\documentclass[11pt]{article}
\usepackage[margin=1in]{geometry}
\usepackage{amsfonts,amsmath,amssymb,bm}
\usepackage{color,cases,float}
\usepackage{graphicx}
\usepackage{hyperref}
\usepackage{psfrag}
\usepackage{eucal}
\usepackage{overpic}
\usepackage{subfigure}
\usepackage{url}
\newcommand{\E}{\mathrm{E}}

\newcommand{\R}{\mathrm{Re}}

\def\tat#1{{\color{black}#1}}

\usepackage{lipsum}

\def\la{\langle}
\def\ra{\rangle}

\def\R{\mathbb{R}}

\def\e{\epsilon}
\def\Z{\mathbb{Z}}
\def\zb{\mathbf{z}}
\def\xb{\mathbf{x}}

\def\mb{\mathbf{m}}
\def\bxi{\boldsymbol{\xi}}

\def\ab{\mathbf{a}}
\def\ba{\boldsymbol{a}}

\def\yb{\mathbf{y}}
\def\bmu{\boldsymbol{\mu}}

\def\Mb{\boldsymbol{M}}
\def\Qb{\boldsymbol{Q}}

\def\Gb{\boldsymbol{G}}

\def\Rb{\boldsymbol{R}}
\def\Ab{\boldsymbol{A}}
\def\Pb{\boldsymbol{P}}

\def\s{\sigma}
\def\g{\gamma}

\def\gb{\mathbf{g}}
\def\R{\mathbb{R}}

\def\Z{\mathbb{Z}}
\def\zb{\mathbf{z}}
\def\xb{\mathbf{x}}

\def\ab{\boldsymbol{a}}
\def\mb{\boldsymbol{m}}
\def\by{\boldsymbol{y}}

\def\bx{\boldsymbol{x}}
\def\Proj{\mbox{Proj}}

\def\r{~}

\definecolor{darkblue}{rgb}{0.0,0.0,0.6}
\hypersetup{colorlinks,breaklinks,
	linkcolor=darkblue,urlcolor=darkblue,
	anchorcolor=darkblue,citecolor=darkblue}
\newtheorem{assumption}{Assumption}
\newtheorem{definition}{Definition}
\newtheorem{lem}{Lemma}
\newtheorem{rem}{Remark}
\newenvironment{proof}{\paragraph{Proof:}}{\hfill$\square$}
\newtheorem{prop}{Proposition}
\newtheorem{theorem}{Theorem}
\newtheorem{cor}{Corollary}
\newtheorem{example}{Example}

\usepackage{algorithm}
\usepackage[noend]{algpseudocode}

\begin{document}

\title{\LARGE
 	Bandit Learning in Convex Non-Strictly Monotone Games
}
\author{Tatiana Tatarenko\thanks{T. Tatarenko is with the Department of Control Theory and Robotics, TU Darmstadt, Germany (e-mail: tatiana.tatarenko@rmr.tu-darmstadt.de)} \and Maryam Kamgarpour\thanks{M. Kamgarpour  is with the EPFL School of Engineering,  Sycamore Lab (e-mail: maryam.kamgarpour@epfl.ch).}\footnotemark[1]
}

\maketitle

\begin{abstract}
We address learning Nash equilibria in convex games under the payoff information setting.  We consider the case in which the game pseudo-gradient is monotone but not necessarily strictly monotone. This relaxation of strict monotonicity enables application of learning algorithms to a larger class of games, such as, for example, a zero-sum game with a merely convex-concave cost function. We derive an algorithm whose iterates provably converge to the least-norm Nash equilibrium in this setting. {From the perspective of a single player using the proposed algorithm,  we view the game as an instance of  online optimization}. Through this lens, we quantify the regret rate of the algorithm and provide an approach to choose the algorithm's parameters to minimize the regret rate.
\end{abstract}

\section{Introduction}

Game theory is a powerful framework to optimize and learn the behavior of multiple interacting agents referred to as players. Such multi-agent problems arise in  application domains including  traffic networks, internet, auctions, and adversarial learning. In several applications, each  player might not know the functional form of her objective. For example, the travel times of different routes in a traffic network \cite{netrout1}, price functions in an electricity market \cite{windfarm,elmark,BasharSG}, or outcomes of an  electricity market auction \cite{karaca2020no} are unknown a priori and depend  on unknown actions of other players.  {By playing the game, a player receives cost
function evaluations at joint played actions, namely, she obtains \emph{payoff information}}\footnote{ {In the optimizaiton and learning community, this information setting is  referred to as zeroth-order  or bandit information/oracle.}}. Our goal is to design a  payoff-based algorithm to learn Nash equilibria in continuous-action merely monotone games.

\subsection{Literature review}\label{sec:intrA}
The payoff-based information setting has been well-explored in \emph{finite action} games in control \cite{vehicle,marden2012revisiting,netrout1} and learning communities \cite{cesa2006prediction}. Payoff-based learning of Nash equilibria in continuous action convex games has been addressed in \cite{bravo2018bandit,tat_kam_TAC, Pang21}. Whereas \cite{Pang21} used two-point function evaluations to estimate gradients, \cite{bravo2018bandit,tat_kam_TAC} leveraged the idea of estimating  the gradient of a player's cost function with respect to her decision variable using only one-point function evaluations, namely, the payoff at a  joint action of all players. This idea is based on \emph{randomized sampling technique} employed  in zeroth-order and stochastic optimization \cite{flaxman2005online, NesterovSpokoiny}. To establish convergence of the iterates to the Nash equilibrium, the \emph{strict monotonicity} of the \emph{pseudo-gradient} of the game was used in \cite{bravo2018bandit,tat_kam_TAC}.

The assumption of strict monotonicity  of the pseudo-gradient rules out several classes of games,  including the class of zero-sum games with convex-concave objective and non-zero sum merely monotone games. Furthermore, games with strictly monotone pseudo-gradient but with an addition of convex coupling constraints \cite{facchinei2007generalized,tat_kam_TAC} have at best a merely monotone  pseudo-gradient.

Leveraging the equivalence of a Nash equilibrium in convex games with the solution set of a so-called  \emph{variational inequality} problem, one can use approaches for solving merely  monotone variational inequalities, such as those based on double time-scale Tikhonov regularization (double time-scale method here refers to an algorithm with an inner and outer loop) or on an extra-gradient  step \cite[Chapter 12]{FaccPang2}. These approaches have been generalized to the stochastic first-order feedback setting, that is to the setting, where stochastic samples of the pseudo-gradient are available \cite{hsieh2019convergence,Lan2022}. A single time-scale approach based on Tikhonov regularization was proposed in \cite{Koshal} with stochastic first-order feedback.  {Note that even though based on the payoff information one can form a stochastic estimate of the game pseudo-gradient at a given point (using randomized sampling), such an estimate would not satisfy the  {bounded variance assumptions} of stochastic first-order feedback and, hence, the analysis does not carry over to the payoff-based setting.} In addition, in all the above work,  the convergence rate of the iterates has  not been addressed. Indeed, past work has used  strong monotonicity of the {pseudo-gradient} to establish convergence rate of the iterates \cite{bravo2018bandit,tat_kam_TAC}.  {Observe that merely monotone games include the subclass of convex optimization. Without additional assumptions, we cannot derive a \emph{rate of iterates' convergence} from that of the function value convergence in the convex optimization setting}. For merely monotone games, an approach to quantify convergence rate is based on the  so-called gap  (or error function) of a variational inequality \cite{cai2022finite,hsieh2019convergence,Nesterov2007} or, alternatively, one may consider the regret rate of a learning algorithm as will be discussed below.

 {From a learning perspective, a fundamental problem is characterizing the class of games for which a so-called \emph{no-regret algorithm}, when employed simultaneously by all players, converges to a Nash equilibrium. In continuous-action games, it has been shown that the mirror-descent class of no-regret algorithms ensure iterates' convergence in potential \cite{heliou2017learning} and strictly monotone games \cite{bravo2018bandit}, even with payoff-based information. However, the same class of algorithms can have divergent iterates in merely monotone games \cite{mertikopoulos2018cycles} even under perfect feedback, that is knowledge of exact pseudo-gradients. More recently, \cite{cai2022finite} proposed a class of no-regret algorithms in merely monotone games with provable convergence to a Nash equilibrium under perfect feedback information. The work in \cite{hsieh2022no} optimized the regret rate for algorithms that have guaranteed convergence to a Nash equilibrium in merely monotone games,  under stochastic first-order feedback. The design of a no-regret algorithms with provable convergence to a Nash equilibrium under  \emph{payoff-based information} in \emph{merely monotone} games was not addressed to our knowledge.}

To our knowledge, the first work addressing payoff-based learning in merely monotone games was our conference paper \cite{tatarenko2019learning}. The approach was inspired by the single time-scale Tikhonov regularization of \cite{Koshal}. More recently, the work \cite{gao2022bandit} addressed learning Nash equilibria  in the  class of games with merely variationally stable equilibria, a class including merely monotone games. In contrast to \cite{tatarenko2019learning}, the algorithm in \cite{gao2022bandit} required memory and focused on the games with interior Nash equilibria. Furthermore,  the no-regret property of the algorithm was not established, and no convergence rate was provided.  This work extends our conference submission to propose an algorithm with  feasible queried actions,   proves that the proposed algorithm is no-regret and provides the setting for the parameters to optimize the regret rate of our algorithm.

To address feasibility of actions, it becomes necessary to modify \cite{tatarenko2019learning} so as to project the algorithm iterates  onto a time-varying shrunk feasible set, while the queried actions are projected onto the original feasible set.  {This modification is needed to ensure the variance of the payoff-based gradient estimator has a suitable order, see Proposition 1.} By properly adjusting this shrinkage parameter along with the Tikhonov regularization one, the randomized sampling distribution, and the stepsize, we ensure convergence of the algorithm iterates, see  Theorem~\ref{th:main}. The addition of this extra shrinkage parameter necessitates new analysis --  {Propositions~\ref{th:Tikhonov} and \ref{prop:zt} extend respectively  analogous results  in \cite[Theorem 12.2.3]{FaccPang1} and \cite[Lemma 3]{koshal2013regularized}, from a single regularied procedure  to a doubly regularized one.}  From the perspective of a single-player, we show that the algorithm is no-regret and quantify the regret rate of the algorithm. In doing so, we provide an approach in setting the algorithm parameters, see Theorem~\ref{th:no-regret_regular} and Corollary~\ref{cor:regret_rate}.

\subsection{Contribution}
In summary, our contributions are as follows. 

\begin{itemize}
\item We develop, to our knowledge, the first convergent memory-free online  algorithm in a convex game with merely monotone pseudo-gradient under the payoff information setting; 
\item We show that our algorithm satisfies the no-regret property when viewed from the perspective of a single player, and regardless of other players' choice of algorithm;
\item We quantify the regret rate of our algorithm. In doing so, we provide an approach to set the algorithm parameters to ensure fast learning quantified in terms of regret rate for a given single player, and convergence to a Nash equilibrium, if employed by all players. 

\end{itemize}

The rest of the manuscript is organized as follows. In Section \ref{sec:problem} we formulate the problem and provide the background on a convex game with a monotone pseudo-gradient. In Section \ref{sec:analysis} we present the algorithm. In Section \ref{sec:main} we derive conditions on the algorithm parameters to ensure convergence of the algorithm iterates to a least-norm Nash equilibrium. In Section \ref{sec:regret} we formulate the online convex program solved by each player and quantify the regret rate of the algorithm. In Section \ref{sec:sim} we evaluate the performance of the proposed algorithm empirically by simulating the algorithm for three  merely monotone games. We conclude in Section \ref{sec:conclusions}. 
 
\subsection{Notations} The set $\{1,\ldots,N\}$ is denoted by $[N]$. Boldface is used to distinguish between vectors in a multi-dimensional space and scalars.
Given $N$ vectors $\bx^i\in\R^d$, $i\in[N]$, $(\bx^i)_{i=1}^{N}$ denotes $({\bx^1}^{\top}, \ldots, {\bx^N}^{\top})^{\top} \in \R^{Nd}$ and $\bx^{-i}:=({\bx^1}^{\top}, \ldots, {\bx^{i-1}}^{\top}, {\bx^{i+1}}^{\top}, \ldots, {\bx^N}^{\top}) \in \R^{(N-1)d}$.
$\R^d_{+}$ denotes vectors from $\R^d$ with non-negative coordinates and   $\Z_{+}$ denotes non-negative whole numbers. The identity matrix in $\R^d$ is denoted by $\mathbf{I}_d$. The identity operator is denoted by $I(\cdot)$. 
The standard inner product on $\R^d$ is denoted by $\la\cdot,\cdot\ra$: $\R^d \times \R^d \to \R$, with associated norm $\|\bx\|:=\sqrt{\la\bx, \bx\ra}$.  {We say the function $f(x): \R\to\R$ is $O(\gb(x))$ as $x\to a$, $f(x) \leq O(g(x))$ as $x\to a$, if $\lim_{x\to a}\frac{|f(x)|}{|g(x)|}\le K$ for some positive constant $K$. } The boundary of a set $Y \subset \R^d$ is denoted by $\partial Y$. For $x \in \R^d$ and a convex closed set $Y \subset \R^d$, $\Proj_{Y}x$ denotes the projection of $x$ onto $Y$. The set of points within $\rho$ distance of the boundary is denoted by $(1-\rho)Y :=\{\bx\in Y: \mbox{dist}(\bx,\partial Y)\ge \rho\}$ with $0< \rho<1$. The indicator function of the set $Y$ is denoted by $\boldsymbol{1}_{\by \in Y}$, which is equal to 1 if $y \in Y$ and equal to $0$ otherwise. The interior of the set $Y$ is denoted by $\mbox{int}(Y)$. 
Expectation of a random variable $\xi$ is denoted by $\E\xi$, whereas its the conditional expectation in respect to some $\sigma$-algebra $\EuScript F$ is denoted by $\E\{\xi\,|\EuScript F\}$.

\section{Background and problem setup}\label{sec:problem}
Consider a game $\Gamma (N, \{\Ab^i\}, \{J^i\})$ with $N$ players, $i$-th player's action set as $\Ab^i\subseteq \R^d$, $i\in[N]$, and her cost (objective) function $J^i:\Ab\to\R$, where $\Ab = \Ab^1\times\ldots\times \Ab^N$ denotes the set of joint actions. 

\begin{definition}\label{def:NE}
An action $\ab^*\in\Ab$ is called a \emph{Nash equilibrium} if  $J^i(\ab^{i*},\ab^{-i*})\le J^i(\ab^{i},\ab^{-i*})$, $\forall i\in[N]$, $\ab^i\in \Ab^i$.
\end{definition}

 
\begin{definition}\label{def:NE}
The \emph{payoff-based information/feedback setting} corresponds to a setting in which each player  knows neither the number of players nor the payoff function of any player (including herself). After each round of play, each player sees its own payoff, that is player $i$ receives $J^i(\ab^i, \ab^{-i})$, where $\ab^i$, $\ab^{-i}$ is hers and others'  chosen actions, respectively. But she  sees neither the choices of other players nor the resulting payoffs.  
\end{definition}


\subsection{Nash Equilibria in convex monotone games}
\begin{assumption}\label{assum:convex}
	$\forall i\in[N]$ the set $\Ab^i$ is convex and compact with a nonempty interior $\mbox{int}(\Ab^i)$, the cost function $J^i(\ab^i, \ab^{-i})$ is defined on $\R^{Nd}$, is continuously differentiable in $\ab$, and convex in $\ab^i$ for  fixed $\ab^{-i}$ on $\R^{Nd}$.
\end{assumption}
Several applications  satisfy these assumptions. For example, in mixed strategy extension of finite action games, the action sets are simplexes and hence, are compact and convex\footnote{Nonempty interior can be guaranteed by reducing the dimension $d$ of the action set and expressing the last coordinate $a^i_d$ as $a^i_d = 1 - \sum_{k\ne d}a^i_k$, see Example 1.}, whereas the cost function of each player is linear in her decision variable. Further examples include Cournot games \cite{FaccPang1}, games arising in wireless communication networks \cite{scutari2012monotone}, electricity markets \cite{chen2014autonomous}, electric vehicle charging \cite{ma2011decentralized}, and route choices in transportation networks \cite{dafermos1980traffic}.

Our approach in designing a payoff-based algorithm lies in connecting Nash equilibria of the game $\Gamma $ with the solution set of a variational inequality problem derived from the game. To state this connection, we first define the game pseudo-gradient and then discuss the variational inequality problem.
\begin{definition}
\label{def:pseudo}
	The  \emph{game pseudo-gradient} $\Mb:\R^{Nd}\to\R^{Nd}$ is defined as
	\begin{align*}
		\nonumber
		\Mb(\ab) &= (\nabla_{\ab^i} J^i(\ab^i, \ab^{-i}))_{i=1}^N=(\Mb^i(\ab))_{i=1}^N, \mbox{where }\\ \nonumber
		\Mb^i(\ab) &= (M^i_1(\ab), \ldots, M^i_d(\ab))^{\top}, \; M^i_{k}(\ab)= \frac{\partial J^i(\ab)}{\partial a^i_k},\; \ab\in\Ab,\\ 
&\text{for all }  i\in[N], \; k\in[d].
	\end{align*}
\end{definition}
In words, the game pseudo-gradient vector stacks the gradient of each player's cost function with respect to her own decision variable. 
If Assumption~\ref{assum:convex} holds, then $\ab^*$ is a Nash equilibrium in $\Gamma (N, \{\Ab^i\}, \{J^i\})$ if and only if $\ab^*$ belongs to the set of solutions to the so-called \emph{variational inequality problem}, denoted by $VI(\Mb, \Ab)$~\cite{FaccPang1}. This solution set is defined as 
\begin{align}\label{eq:VI}
SOL(\Ab,\Mb):= \{\ab^* \,| \, \la\Mb(\ab^*), \ab - \ab^*\ra\ge 0, \quad \forall \ab\in\Ab\}.
\end{align}


\begin{definition}\label{def:mm}
	A mapping $\Mb:\R^n\to\R^n$ is \emph{monotone} over $Y\subseteq\R^n$ if $\la\Mb(\by)-\Mb( \tilde{\by}),\by- \tilde{\by}\ra\ge 0$ for every $\by, \tilde{\by}\in Y$. It is \emph{strictly monotone} if  $\forall \by \neq \tilde{\by}$ the above inequality is strict and \emph{strongly monotone} if  there exists $\tau > 0$ such that $\la\Mb(\by)-\Mb( \tilde{\by}),\by- \tilde{\by}\ra\ge \tau \| \by - \tilde{\by} \|^2$. If $\Mb$ is monotone but not strictly monotone, we refer to it as \emph{merely monotone}. Furthermore, we refer to a game with (strongly/strictly/merely) monotone pseudo-gradient as a (strongly/strictly/merely) monotone game. 
\end{definition}
As an example, given a linear mapping $\Mb(\ab) = H \ab + h$, $H \in \R^{Nd\times Nd}$  then $\Mb$ is strongly monotone if the symmetric part of $H$, $H^s :=\frac{1}{2}(H+H^T)$, is positive definite, and merely monotone if $H^s$ is positive semi-definite. 
\begin{assumption}\label{assum:CG_grad}
The  mapping $\Mb$ is \emph{merely monotone on $\Ab$}.
\end{assumption}
The monotonicity condition is satisfied in a large class of games  \cite{Rosen1965,chen2014autonomous,ma2011decentralized,scutari2012monotone}.

\subsection{Examples}\label{sec:example}
The  examples below illustrate zero-sum and non zero-sum merely monotone games, the non-uniqueness of Nash equilibria in this class of games, and that player-wise (strict/strong) convexity is not sufficient for (strict/strong) monotonicity. 

\begin{example}\label{example:penny}
	Consider the zero-sum game of matching pennies. Each player has two actions $\{0,1\}$. The cost for player 1 corresponding to actions $i,j \in \{0,1\}$ of player one and two, respectively, is denoted by $c_{i,j}$ and is given by $c_{0,0} = c_{1,1} = 1$ and $c_{0,1} = c_{1,0} = -1$. By relaxing the action sets to $\Ab^i = [0,1] \subset \R$, where, $a^i \in \Ab^i$ denotes the probability of playing action $0 \in \{0,1\}$, we obtain the mixed strategy extension of the game. The payoff of the mixed strategy game is the expected payoff of the finite action game: $J^1(a^1, a^2) = a^1 c_{0,0} a^2 +a^1 c_{0,1} (1-a^2) + (1-a^1) c_{1,0} a^2 + (1-a^1) c_{1,1} (1-a^2)=4a^1a^2 -2(a^1+a^2) +1$, $J^2 = -J^1$.  Clearly, the game satisfies Assumptions~\ref{assum:convex}. The pseudo-gradient is given by $\Mb(a^1, a^2) = (4a^2-2, -4a^1+2)^T$. It is merely monotone, and hence the game satisfies Assumption~\ref{assum:CG_grad}. The unique   Nash equilibrium is at $(1/2, 1/2)$. Considering a modified game with the same payoff but  restricting the action sets to $\Ab^i = [1/2,1]$, $i=1,2$,  it can be verified that the Nash equilibria are $(1/2, a^2)$ for any $a^2 \geq 1/2$. The least-norm Nash equilibrium is  $(1/2, 1/2)$. We will be referring to the cases with action sets $A^i=[0,1]$,  $A^i=[1/2,1]$ for $i=1,2$, as Examples 1a and 1b, respectively. 
\end{example}

\begin{example}\label{example:nonunique}
	Consider a non zero-sum game with $\Ab^1 = \Ab^2 = [-1, 1] \subset \R$ and $J^1(a^1,a^2) = a^1a^2 + \frac{1}{2}(a^1)^2$ and $J^2(a^1,a^2) = a^1a^2 + \frac{1}{2}(a^2)^2$. The game pseudo-gradient is given by $\Mb(a^1, a^2) = (a^1+a^2, a^1+a^2)^T$, which is merely monotone.  {Observe however that the cost of each player is strongly convex in her decision variables.}  The game satisfies Assumptions~\ref{assum:convex},~\ref{assum:CG_grad}.  This game has a continuum of Nash equilibria, namely, $a^1 = -a^2$ is a Nash equilibrium, $\forall a^1 \in \Ab^1$. The least-norm Nash equilibrium is $a^1 = a^2 = 0$.
\end{example}

\section{Proposed Learning Algorithm}\label{sec:analysis}

The idea of the proposed algorithm is that a given player $i$  estimates the gradient of its cost function $J^i$  at a played action $\ab^i$ using her payoff information $J^i(\ab^i, \ab^{-i})$. Then, she performs a doubly regularized approximate gradient-descent, where the first regularization addresses mere monotonicity, and the second addresses feasibility of the iterates. We  formalize this approach in Section~\ref{sec:reg_alg}. In Section \ref{sec:est_prop} we derive the desired properties of the gradient estimation procedure.

\subsection{Decoupled payoff-based regularized learning}\label{sec:reg_alg}

\begin{algorithm}[t!]
	\caption{Payoff-based learning of in merely monotone games}\label{alg:algorithm1}
	\begin{algorithmic}[1]
		\Require Action set $\Ab^i \subset \R^d$, the sequences $\{\g_t\}, \{\sigma_t\}, \{\rho_t\}, \{\e_t\}$, initial iterate $\bmu^i(0)$.
		\For {$t = 0,1, \ldots$}
		\State  Sample $\bxi^i(t)$ according to the probability density function defined in Equation~\eqref{eq:density}.
		\State Play the feasible action $\ab^i(t) = \Proj_{\Ab^i}[\bxi^i(t)]$.
		
		\State Receive $ J^i(t) = J^i(\ab^1(t),\ldots,\ab^N(t))$.
		\State Estimate  local gradient of $J^i$ as per Equation \eqref{eq:est_Gd}.
		\State Set $\bmu^i(t+1)$ as per per Equation \eqref{eq:rppg_i}.
		
		\EndFor
		
	\end{algorithmic}
	\label{alg:payoff}
\end{algorithm}

Denote by $\mb^i$ player $i$'s estimate of the term $\Mb^i$ in the pseudo-gradient of the game (see Definition \ref{def:pseudo}). The proposed procedure  to update player $i$'s \emph{iterate}, denoted by $\bmu^i$, is as:
\begin{align}
	\label{eq:rppg_i}
	\bmu^i(t+1)=\Proj_{(1-\rho_t)\Ab^i}[\bmu^i(t)-\gamma_t(\mb^i(t)+\e_t\bmu^i(t))],
\end{align}
where $\bmu^i(0)\in \R^{Nd}$ is an arbitrary finite value, $\g_t$ is the stepsize, and  $\e_t$, $\rho_t$ are the two regularization parameters. The parameter $\rho_t\in[0,1)$ is introduced to control the feasibility of the played actions, whereas $\e_t$ introduces a regularization to address the merely monotone pseudo-gradient. The term $\mb^i(t)$ is obtained using the payoff-based feedback as described below.

Given $\bmu^i(t)$, let  player $i$ sample the random vector $\bxi^i(t)$ according to the multidimensional normal distribution $\EuScript N(\bmu^i(t)=(\mu^i_1(t),\ldots,\mu^i_{d}(t))^{\top},\sigma_t)$ with the density function
\begin{align}\label{eq:density}
	p^i&(\bx^i;\bmu^i(t),\sigma_{t})= \frac{1}{(\sqrt{2\pi}\sigma_{t})^{d}}\exp\left\{-\sum_{k=1}^{d}\frac{(x^i_k-\mu^i_k(t))^2}{2\sigma^2_{t}}\right\}.
\end{align}
Then, the action agent $i$ plays is $\ab^i(t) = \Proj_{\Ab^i}[\bxi^i(t)]$.  According to the payoff information setting under consideration, the cost value $ J^i(t)$ at the joint action $\ab(t)=(\ab^1(t),\ldots,\ab^N(t))\in \Ab$,  denoted by $J^i(t) :=J^i(\ab(t))$, is revealed to each player $i$.  {Observe that $\bmu^i(t)$ is the the iterate of the algorithm, in  contrast  with the actual played actions $\ab^i(t)$ and it satisfies the dynamics as per Equation~\eqref{eq:rppg_i}}. \tat{Taking into account the procedure~\eqref{eq:rppg_i} to update the iterate $\bmu^i(t)$}, each player $i$ then estimates her local gradient $\Mb^i$ evaluated at the joint iterate $\bmu(t)=(\bmu^i(t))_{i=1}^N$ as: 
\begin{align}\label{eq:est_Gd}
	\mb^i(t) = { J^i(t)}\frac{{\bxi^i(t)} -\bmu^i(t)}{\sigma^2_t}.
\end{align}
Putting these steps together, Algorithm \ref{alg:algorithm1} specifies the payoff-based learning algorithm for each player $i$.
 {\begin{rem}
Player $i$ can have her local choice of the parameters $\sigma^i_t, \rho^i_t, \g^i_t, \e^i_t$. For simplicity, we drop the dependence of these parameters on $i$. However, the analysis can be extended to this case (see, for example, \cite{tat_kam_TAC} a similar extension but for strictly monotone games).
\end{rem}}
\begin{rem}\label{rem:distribution_choice}
Our choice of Gaussian distribution is inspired by stochastic sampling in \cite{Thatha, NesterovSpokoiny}. One can choose to estimate the gradient by sampling the actions from the uniform distribution on a sphere \cite{flaxman2005online,bravo2018bandit}, see Appendix~\ref{appendix:uniform}  for details. 
In either case, the bias and variance of the estimates will be of similar order. Importantly, it is  the order of these terms that will determine the convergence of the algorithm.
\end{rem}

\subsection{Gradient estimator properties}\label{sec:est_prop}
The following assumption is needed to bound the bias and variance of the gradient estimator defined in Equation \eqref{eq:est_Gd}.

\begin{assumption}
	\label{assum:Lipschitz}
	\begin{enumerate}
		\item The pseudo-gradient $\Mb$ is continuously differentiable over $\R^{Nd}$;\label{itm:lip}
		\item For each function $J^i$ the following holds: $J^i(\bx) = O(\exp\{\|\bx\|^{\alpha}\})$ as $\|\bx\|\to\infty$, where $\alpha<2$. \label{itm:growth}
	\end{enumerate}
\end{assumption}


\begin{rem}\label{rem:LipschJ}
In addition to  Lipschitz continuity of $\Mb$ employed in literature on learning with payoff information~\cite{bravo2018bandit, tat_kam_TAC}, in part~\ref{itm:lip}) of the assumption above we  assume differentiability of $\Mb$  to bound the error arising due to the projection onto the shrunk set $(1-\rho_t) \Ab^i$, see  proof of Proposition~\ref{prop:sample_grad}.
Part~\ref{itm:growth} of Assumption~\ref{assum:Lipschitz}  addresses the unbounded Gaussian distribution employed in querying the actions. In particular, this assumption is used to upper bound the moments of the functions $J^i$, $i\in[N],$ and those of their gradients in the proof of the main results, Theorems~\ref{th:main},~\ref{th:no-regret_regular}.
\end{rem}

The next lemma establishes that by using  $\mb^i$, defined in \eqref{eq:est_Gd}, each agent can obtain an estimate of its gradient $\Mb^i$, with a bias term $\Qb^i$ and a zero-mean noise term $\Rb^i$. 
Denote $\EuScript F_t$ as the $\sigma$-algebra generated by the random variables $\{\bmu(k),\bxi(k)\}_{k\le t}$. 
\begin{prop}\label{prop:sample_grad}
	Consider the game $\Gamma (N, \{\Ab^i\}, \{J^i\})$ for which Assumptions~\ref{assum:convex} and ~\ref{assum:Lipschitz} hold. Choose $\sigma_t, \rho_t$ such that $\lim_{t\to\infty}\sigma_t = 0, \lim_{t\to\infty}\rho_t = 0, \lim_{t\to\infty}\frac{\rho_t}{\s_t}=\infty $. Then the gradient estimate $\mb^i(t)$ introduced in Equation \ref{eq:est_Gd} can be decomposed as follows:
	\begin{align}
		\label{eq:map_est}
		{\mb^i}(t) = \Mb^i(\bmu(t)) + \Qb^i(t) + \Rb^i(t),
	\end{align}
	where the stochastic terms $\Qb^i(t)$ and $\Rb^i(t)$ satisfy
	\begin{enumerate}
		\itemsep-.3em
		\item $\E\{\Rb^i(t) | \EuScript F_t\} =0$,  \label{itm:R_mean}
		\item $ \E\{\|\Rb^i(t)\|^2 | \EuScript F_t\} = O\left(\frac{1}{\sigma_t^2}\right)$,  \label{itm:R_square}
		\item $\E\{\|\Qb^i(t)\||\EuScript F_t\} = O(\sigma_t)$.  \label{itm:Q_norm}
	\end{enumerate}
\end{prop}

\begin{proof}
Let us define the smoothed cost function for player $i$ as $\tilde{J}^i: \R^{Nd} \rightarrow \R$, namely
\begin{align}
	\label{eq:mixedJ}
	\tilde{J}^i &(\bmu^1,\ldots,\bmu^N)= \int_{\mathbb R^{Nd}}J^i(\bx)p( \bx; \bmu, \sigma)d\bx,  
\end{align}
where $p( \bx; \bmu, \sigma) = \prod_{i=1}^{N}p( \bx^i; \bmu^i, \sigma)$ is the joint probability density function.
In other words, $\tilde{J}^i$ is the cost in mixed strategies, given the actions are distributed according to the density function in  \eqref{eq:density}.  \tat{Denote  by $\tilde{\Mb}^i(\bmu):\R^{Nd} \rightarrow \R^{Nd}$ the terms in the pseudo-gradient (see Definition \ref{def:pseudo}) of the game with smoothed costs  at $\bmu=(\bmu^1,\ldots,\bmu^N)$}, i.e.,
\begin{align}\label{eq:smoothed_main}
	\tilde{\Mb}^i(\bmu) =\frac{\partial \tilde{J}^i(\bmu^1,\ldots,\bmu^N) }{\partial \bmu^i},  \; i \in[N].
\end{align}
With simple addition and subtraction, one can verify that the terms in Equation~\eqref{eq:map_est} can be defined as:
\begin{align}
	\label{eq:Qi0}
	&\Qb^i(t) :=\underbrace{\tilde{\Mb}^i (\bmu(t)) -\Mb^i(\bmu(t))}_{\text{term 1}} \cr
	&\qquad+ \underbrace{(J^i(\ab(t))-J^i(\bxi(t)))\frac{\bxi^i(t) -\bmu^i(t)}{\sigma^2_t}}_{\text{term 2 = $\Pb^i(t)$}},\\
	\label{eq:Ri_main}
	&\Rb^i(t) := J^i(\bxi(t))\frac{\bxi^i(t) -\bmu^i(t)}{\sigma^2_t} - \tilde{\Mb}^i (\bmu(t)).
\end{align}
Having the definitions above in place, parts 1) and 2) of the proposition above are an easy extension of the result in Lemma~1 and estimation~(23) in~\cite{tatarenko2019learning} - this extension is provided for completeness in Appendix \ref{app:sample_grad}. Term 1 in the definition~\eqref{eq:Qi0} of $\Qb^i(t)$ is estimated also in~\cite{tatarenko2019learning} as: $\|\tilde{\Mb}^i (\bmu(t)) -\Mb^i(\bmu(t))\| = O(\sigma_t)$ (see the estimation Equation (20) therein).
Thus, we focus on proving part 3) of the proposition by bounding term~2 in~\eqref{eq:Qi0}, denoted as $\Pb^i(t)$, namely $\Pb^i(t) = (J^i(\ab(t))-J^i(\bxi(t)))\frac{\bxi^i(t) -\bmu^i(t)}{\sigma^2_t}$. This  term is new compared to \cite{tat_kam_TAC, tatarenko2019learning} and arises due to the projection of the sampled vectors $\bxi^i(t)$ onto the feasible sets.  
  
 {
 	\begin{align}\label{eq:projterm0}
  &\E\left\{\Pb^i(t)\;|\EuScript F_t\right\}\cr
 		&=\E\left\{|J^i(\ab(t))-J^i(\bxi(t))|\frac{\|\bxi^i(t) -\bmu^i(t)\|}{\sigma^2_t}\;|\EuScript F_t\right\}
   \cr
&\le (\E\left\{|J^i(\ab(t))-J^i(\bxi(t))|^2\right\})^{\frac{1}{2}}\left(\E\left\{\frac{\|\bxi(t)-\bmu(t)\|^2}{\sigma^4_t}\right\}\right)^{\frac{1}{2}}\cr
 		& =\frac{\sqrt{Nd}}{\sigma_t}(\E\left\{|J^i(\ab(t))-J^i(\bxi(t))|^2\right\})^{\frac{1}{2}},
 	\end{align}
where we used the H\"older's inequality ($\E|XY|\le (\E X^2)^{1/2}(\E Y^2)^{1/2})$ and the notation $\E\{\cdot\} = \E\{\cdot |\bxi(t)\sim \EuScript(\bmu(t),\sigma_t)\}$.
 	Next, we estimate the term $\E\left\{|J^i(\ab(t))-J^i(\bxi(t))|^2\right\}$ as follows. As $\ab(t) = \Proj_{\Ab}\bxi(t)$, we conclude that almost surely:
 	\begin{align}\label{eq:projterm}
 		&\E\left\{|J^i(\ab(t))-J^i(\bxi(t))|^2\right\} \cr
 		&=\int_{\R^{Nd}\setminus \Ab}|J^i(\Proj_{\Ab}\bx)-J^i(\bx)|^2p(\bx;\bmu(t),\sigma_t)d\bx\cr
 		&\le 2\int_{\R^{Nd}\setminus \Ab}(J^i(\Proj_{\Ab}\bx))^2p(\bx;\bmu(t),\sigma_t)d\bx\cr
   &\,+ 2\int_{\R^{Nd}\setminus \Ab}(J^i(\bx))^2p(\bx;\bmu(t),\sigma_t)d\bx\cr
 		&\le 2K\Pr\{\bxi(t)\in\R^{Nd}\setminus \Ab\}) \cr
   &\,+ 2\int_{\R^{Nd}\setminus \Ab}(J^i(\bx))^2p(\bx;\bmu(t),\sigma_t)d\bx,
 	\end{align}
 	where the first equality above is due to the fact that $\Proj_{\Ab}\bx=\bx$ for any $\bx\in \Ab$. The first inequality was obtained by taking into account that $|J^i(\Proj_{\Ab}\bx)-J^i(\bx)|^2\le 2 (J^i(\Proj_{\Ab}\bx))^2 + 2(J^i(\bx))^2$, whereas the last one is due to the inequality $(J^i(\Proj_{\Ab}\bx))^2\le K$ for some constant $K$ and any $\bx$ (since $\Ab$ is compact and $J^i$ is continuous).}
 
  {	Thus, let us estimate $\Pr\{\bxi(t)\in\R^{Nd}\setminus \Ab\}$ above. 
 	 	Let $\EuScript O_{\rho_t}(\xb) = \{\boldsymbol y\in\R^{Nd}| \|\boldsymbol y-\xb\|^2<{\rho}_t^2\}$ denote the $\rho_t$-neighborhood of the point $\xb\in\Ab$.  Then, taking into account the fact that $\EuScript O_{\rho_t}(\bmu(t))$ is contained in $\Ab$ and $\rho_t<1$, we obtain that $\forall t$ and any $\sigma$ such that $\sigma>\sigma_t$:
{\allowdisplaybreaks
\begin{align}\label{eq:probterm}
 		\Pr&\{\bxi(t)\in\R^{Nd}\setminus \Ab\}\le\Pr\{\bxi(t)\in\R^{Nd}\setminus \EuScript O_{\rho_t}(\bmu(t))\}\cr 
 		& = \int_{\boldsymbol y\notin \EuScript O_{\rho_t}(\bmu(t))}\frac{1}{(2\pi)^{Nd/2}\sigma_t^{Nd}}
 		\exp\left\{-\frac{\|\boldsymbol y-\bmu(t)\|^2}{2\sigma_t^2}\right\}d\boldsymbol y \cr
 		&=\int_{\boldsymbol y\notin \EuScript O_{\rho_t}(\bmu(t))}\exp\left\{-\|\boldsymbol y-\bmu(t)\|^2\left(\frac{1}{2\sigma_t^2} - \frac{1}{2\sigma^2}\right)\right\}\cr
 		&\qquad\qquad\times\frac{\sigma^{Nd}}{\sigma_t^{Nd}}\frac{1}{(2\pi)^{{Nd}/2}\sigma^{Nd}}\exp\left\{-\frac{\|\boldsymbol y-\bmu(t)\|^2}{2\sigma^2}\right\}d\boldsymbol y\cr
 		&\le \exp\left\{-\rho^2_t\left(\frac{1}{2\sigma_t^2} - \frac{1}{2\sigma^2}\right)\right\}\frac{\sigma^{Nd}}{\sigma_t^{Nd}}\cr
 		&\quad\times\int_{\boldsymbol y\notin \EuScript O_{\rho_t}(\bmu(t))}\frac{1}{(2\pi)^{Nd/2}\sigma^{Nd}}\exp\left\{-  \frac{\|\boldsymbol y-\bmu(t)\|^2}{2\sigma^2}\right\}d\boldsymbol y\cr
 		&\le k_1 \frac{e^{-\frac{{\rho}_t^2}{2\sigma_t^2}}}{\sigma_t^{Nd}} = O(f(t)),	\end{align}
}
for some finite $k_1>0$, and $f(t) = \frac{e^{-\frac{{\rho}_t^2}{2\sigma_t^2}}}{\sigma_t^{Nd}}$. The last inequality holds because
\[\int_{\boldsymbol y\notin \EuScript O_{\rho_t}(\bmu(t))}\frac{1}{(2\pi)^{{Nd}/2}\sigma^{Nd}}\exp\left\{-  \frac{\|\boldsymbol y-\bmu(t)\|^2}{2\sigma^2}\right\}d\boldsymbol y\le 1\]
and, thus, due to the diminishing $\rho_t$ there exists $0<k_1<\infty$:
\[\int_{\boldsymbol y\notin \EuScript O_{\rho_t}(\bmu(t))}\frac{e^{ \frac{{\rho}_t^2}{2\sigma^2}}\sigma^{Nd}}{(2\pi)^{N/2}\sigma^{Nd}}\exp\left\{-  \frac{\|\boldsymbol y-\bmu(t)\|^2}{2\sigma^2}\right\}d\boldsymbol y\le k_1.\]
Next, regarding the second term in the last inequality of \eqref{eq:projterm}, taking into account Assumption~\ref{assum:Lipschitz}, part~\ref{itm:growth}), we conclude existence of some constant $C_1$ such that $(J^i(\bx))^2\exp\left\{-\frac{\|\boldsymbol x-\bmu(t)\|^2 }{4\sigma_t^2}\right\}\le C_1$ for any $\bx$. Thus,
 	\begin{align*}
 		&\int_{\R^{Nd}\setminus \Ab}(J^i(\bx))^2p(\bx;\bmu(t),\sigma_t)d\bx\cr
 		& \le C_1\int_{\boldsymbol x\notin \EuScript O_{\rho_t}(\bmu(t))}\frac{1}{(2\pi)^{Nd/2}\sigma_t^{Nd}}
 		\exp\left\{-\frac{\|\boldsymbol x-\bmu(t)\|^2 }{4\sigma_t^2}\right\}d\boldsymbol x.
 	\end{align*}
 	Hence, analogously to the derivation steps of \eqref{eq:probterm}, we obtain (see Appendix \ref{app:sample_grad} in  for more steps detailed out)
 	\begin{align}\label{eq:probterm1}
 		&\int_{\R^{Nd}\setminus \Ab}(J^i(\bx))^2p(\bx;\bmu(t),\sigma_t)d\bx 
   = O\left(f(t)\right),
 	\end{align}
Since $\lim_{t\rightarrow \infty}\frac{\rho_t}{\sigma_t} = \infty$ by the assumption of the proposition, we have that $f(t) = O(\sigma_t)$ and thus, 
term 2, equal to $\Pb^i(t)$ in~\eqref{eq:Qi0}, has the order less than $O(\sigma_t)$ and, thus, $\E\{\|\Qb^i(t)\||\EuScript F_t\} = O(\sigma_t)$.}
\end{proof}



Note that the variance and bias terms in parts 2 and 3 of Proposition~\ref{prop:sample_grad} are  functions of $\sigma_t$ and, thus, can be controlled by appropriately choosing $\sigma_t$. Furthermore, the order of these terms guide us in choosing the stepsize $\g_t$ and regularization terms $\e_t, \rho_t$ to ensure convergence of the algorithm iterates. These choices are detailed in the next section.

\begin{rem}\label{rem:distribution_choice}
 {The proposed procedure in Algorithm~\ref{alg:algorithm1} corresponds to a single time-scale Tikhonov regularization approach. This approach is based on the idea to modify the initial merely monotone pseudo-gradient mapping $\Mb(\cdot)$ to a strongly monotone one $\Mb(\cdot) + \epsilon_t I(\cdot)$ at each iteration, given $\epsilon_t>0$. Then, one iteration of  gradient-descent   is applied player-wise based on the gradient of this modified mapping, and $\epsilon_t$ is adjusted for the next iteration. The intuition is that for a fixed $\epsilon$, the gradient descent approach has convergence due to strong monotonicity
\footnote{Gradient-descent  without any regularization term generally diverges in the case of merely monotone variational inequality problems~\cite{GRAMMATICO2018186}.}. The challenge here is to balance the rate of the regularization parameter $\epsilon_t$ tending to zero  with other algorithm's parameters, so as to guarantee the iterates' convergence to some solution of the initial variational inequality formulated with respect to  the pseudo-gradient $\Mb(\cdot)$.  
}
\end{rem}

\section{Convergence to a Nash Equilibrium}\label{sec:main}
We demonstrate that there exists a set of parameters such that the proposed payoff-based  algorithm converges to a least-norm Nash equilibrium in the game $\Gamma(N, \{\Ab^i\}, \{J^i\})$.  

\begin{assumption}\label{assum:parameter2}
	Choose $\g_t$, $\sigma_t$, $\rho_t$, $\e_t$ to satisfy $\lim_{t\to\infty} \g_t = \lim_{t\to\infty} \e_t = \lim_{t\to\infty} \s_t = \lim_{t\to\infty} \rho_t =0$, and 
	\begin{enumerate}
		\itemsep-.3em
		\item $\lim_{t\to\infty}\frac{\rho_t}{\e_t} = 0$,  \label{itm:rho_eps}
		\item $\sum_{t=0}^{\infty}\left(\frac{|\e_t-\e_{t-1}|^2}{\e^{3}_t\g_t}+\frac{|\rho_t-\rho_{t-1}|^2}{\e^{5+2\varepsilon}_t\g_t}\right)<\infty$ \mbox{ for some $\varepsilon>0$}, \label{itm:omega}
		\item $\lim_{t\to\infty}\frac{\rho_t}{\s_t} = \infty$,   \label{itm:rho_sig}\\
		\item $\sum_{t=0}^{\infty}\g_t\e_t = \infty$,  \label{itm:gam_eps}\\
		\item $\sum_{t=0}^{\infty}\frac{\g^2_t}{\s^2_t}<\infty$, \label{itm:gam_sig}
		\item $\sum_{t=0}^{\infty}\g_t\s_t<\infty$.\label{itm:gam_sig2}
	\end{enumerate}
	Moreover, if the least-norm solutions of $VI(\Mb,\Ab)$ are contained in $\mbox{int}(\Ab)$, relax condition~\ref{itm:omega}) above by
	\begin{enumerate}
		\itemsep-.3em
		\item[2)*] $\sum_{t=0}^{\infty}\frac{|\e_t-\e_{t-1}|^2}{\e^{3}_t\g_t}<\infty.$ \label{itm:omega1}
	\end{enumerate}
\end{assumption}

\begin{theorem}\label{th:main}
	Let each of the players in the game $\Gamma(N, \{\Ab^i\}, \{J^i\})$ update their  iterates $(\bmu^i(t))_{i=1}^N$ according to Algorithm \ref{alg:algorithm1}.  Under Assumptions\r\ref{assum:convex}-\ref{assum:parameter2}, given an arbitrary $\bmu(0)$, the  joint iterate $\bmu(t)=(\bmu^1(t),\ldots, \bmu^N(t))$ converges almost surely to a least norm Nash equilibrium $\ab^*$, and the joint action $\ab(t)$ converges in probability to $\ab^*$.
\end{theorem}

We provide the proof of the theorem in Subsection \ref{sec:proof_sketch}. Let us provide intuition for the conditions in Assumption~\ref{assum:parameter2}. Furthermore, we  demonstrate  in Lemma \ref{lem:step_constants} that there is a continuum  choice of parameters that satisfy Assumption~\ref{assum:parameter2}. 

First, we discuss handling of the error due to approximating the pseudo-gradient $\Mb^i$ using the payoff information. Notice that using payoff-based feedback we obtain a stochastic estimate of the gradient with bias and variance terms, as per Proposition~\ref{prop:sample_grad}. Thus, our algorithm can be interpreted as an analogue of a stochastic approximation procedure. Conditions~\ref{itm:gam_eps}), \ref{itm:gam_sig}), \ref{itm:gam_sig2}) of Assumption~\ref{assum:parameter2} are standard in stochastic approximation \cite{borkar2008stochastic}. In particular, $\sum_{t=0}^{\infty}\g_t\e_t = \infty$ ensures the algorithm can make sufficient progress with time. On the other hand, $\sum_{t=0}^{\infty}\frac{\g^2_t}{\s^2_t}<\infty$ and  $\sum_{t=0}^{\infty}\g_t\sigma_t < \infty$ keep the perturbations caused by the stochastic  nature of the algorithm under control, and are a consequence of  the order of variance and bias terms derived in parts 2 and 3 of Proposition \ref{prop:sample_grad}.

Second, condition~\ref{itm:rho_sig}) on the slower growth of $\rho_t$ relative to $\sigma_t$ arises due to requiring feasibility of actions during the learning and is used to bound term 2 defined in Equation~\eqref{eq:Qi0} arising from projection of the actions onto the feasible set.  It is in particular a consequence of concentration of the Gaussian distribution around its mean.

{Third, we discuss the conditions arising due to mere monotonicity of the pseudo-gradient. In particular, conditions~\ref{itm:rho_eps}) and \ref{itm:omega}) allow for comparing the single time-scale regularized procedure with that of a double time-scale procedure. In this standard double time-scale approach \cite{FaccPang2}  at each inner step  $VI((1-\rho_{t})\Ab, \Mb + \e_{t} I)$ is solved, while $\e_t$ is decreased in the outer step. Moreover, condition \ref{itm:omega}), or its variant in condition~2)* for the case of the interior least-norm  Nash equilibria, arises due to the bound on the norm of the difference between the solutions of the two variational inequalities $VI((1-\rho_t) \Ab, \Mb + \e_t I)$, and $VI((1-\rho_{t-1})\Ab, \Mb + \e_{t-1} I)$. 
As will be  shown (in the proof of Proposition~\ref{prop:zt} in ~\ref{sec:proof_sketch}), if the least-norm Nash equilibria belong to the interior of the set $\Ab$, the convergent series $\sum_{t=0}^{\infty}\frac{|\e_t-\e_{t-1}|^2}{\e^{3}_t\g_t}$ is sufficient to bound this difference, resulting in condition~2)*. Otherwise, the convergence of the series $\sum_{t=0}^{\infty}\frac{|\rho_t-\rho_{t-1}|^2}{\e^{5+2\varepsilon}_t\g_t}$, which is related to the behavior of projection onto the shrunk sets, needs to be accounted for, resulting  in condition~\ref{itm:omega}).

We now derive an approach to set the parameters so that they satisfy Assumption~\ref{assum:parameter2}. 
\begin{lem}\label{lem:step_constants}
	Consider $\gamma_t=O(\frac{1}{t^g})$, $\sigma_t = O(\frac{1}{t^s})$,  $\rho_t = O(\frac{1}{t^r})$, $\e_t=O(\frac{1}{t^e})$, $0< g,e, s,r < 1$.
	A sufficient condition for satisfying Assumption~\ref{assum:parameter2} is choosing $g, e, s, r$ to satisfy:
	\begin{enumerate}
		\item $g+5e-2r<1$,\label{cons:ger}
		\item $e< r < s, \; g + e \le 1$, \label{cons:ge}
		\item $2g-2s>1$, \; $g+s > 1$. \label{cons:gs}
	\end{enumerate}
	Moreover, if  the least-norm solutions of $VI(\Mb,\Ab)$ are contained in $\mbox{int}(\Ab)$, the first condition above is not needed.
\end{lem}
The proof of the above lemma is based on the the convergence of a $p$-series and is provided in (see Appendix \ref{app:step} therein).

{The above conditions define a set of linear constraints with non-empty interior for the parameters.  As an example of a set of parameters satisfying Lemma~\ref{lem:step_constants}, one can pick $(g,s,r,e)=(0.79, 0.25, 0.23, 0.21 )$  in the case in which the least norm solutions of $VI(\Mb, \Ab)$ are in the interior of $\Ab$} and $(g, s, r, e)=(0.87, 0.33, 0.29, 0.13)$ otherwise (these choices of parameters are based on the regret rate derived in the next section).  {Other feasible parameters can be visualized from the polytope of feasible set,  by fixing some of the parameters and considering the feasible region of the remaining parameters. The feasible regions in $r, e$ plane  corresponding to fixing $g,s$ based on the two vectors above are shown in Figure \ref{fig:feasible}. }

\begin{figure}[!t]
\begin{subfigure}
		\centering
		\includegraphics[width = .8\linewidth]{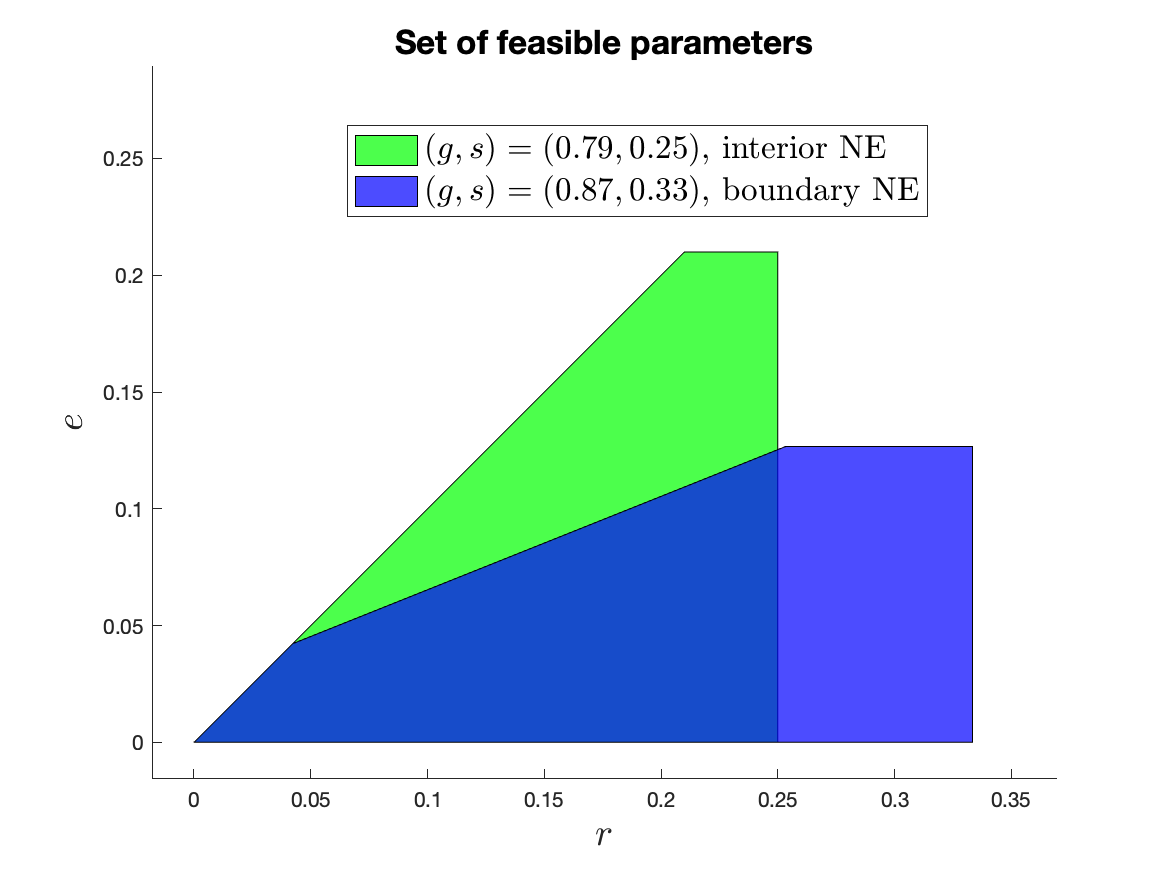}
		\caption{An example of set of $r$,$e$ satisfying Assumption \ref{assum:parameter2}, for fixed $g,s$.}
		\label{fig:feasible}
	\end{subfigure}%
\end{figure}

\subsection{Proof of  Theorem \ref{th:main} }\label{sec:proof_sketch}
\begin{proof}
The proof consists of three steps, summarized here.  

\subsubsection*{Step 1} We first introduce the auxiliary sequence $\by(t) = (\by^1(t), \dots, \by^N(t))\in\R^{Nd}$ as the solution of the variational inequality $VI((1-\rho_t)\Ab,\Mb(\by) + \e_t\by)$, namely
\begin{align}\label{eq:yt}
\by(t)\in SOL((1-\rho_t)\Ab,\Mb(\by) + \e_t\by).
\end{align}
The standard Tikhonov sequence corresponds to the solution of the variational inequality above with the $\rho_t =0$ \cite[Chapter 12]{FaccPang1}. We now establish a similar result to  that of the standard Tikhonov sequence. See Appendix~\ref{app:y_proof} of  for proof.
 \begin{prop}\label{th:Tikhonov}
 Under Assumptions~\ref{assum:CG_grad} and~\ref{assum:Lipschitz},  $\by(t)$ defined in \eqref{eq:yt} exists and is unique for each $t$. Moreover, for $\e_t \to 0$ and $\rho_t\to 0$ and given $\lim_{t\to\infty}\frac{\rho_t}{\e_t} = 0$,  $\by(t)$ is uniformly bounded and converges to the least-norm solution $\ab^*$ of $VI(\Ab,\Mb)$.
 \end{prop}

\subsubsection*{Step 2}  Next, consider  the following regularized procedure
	\begin{align}\label{eq:crmd}
		\zb(t+1)=\Proj_{(1-\rho_t)\Ab}[\zb(t)-\gamma_t(\Mb(\zb(t))+\e_t\zb(t))],
	\end{align}
	where $\zb_0\in\R^{Nd}$ is arbitrary.
Letting  $\bmu(t) = (\bmu^i(t))_{i=1}^N$, observe that the iterates of our algorithm \eqref{eq:rppg_i} can be written compactly as
	\begin{align}
		\label{eq:rppg}
		\bmu(t+1)=\Proj_{(1-\rho_t)\Ab}[\bmu(t)-\gamma_t(\mb(t)+\e_t\bmu(t))].
	\end{align}
Thus,~\eqref{eq:crmd} is analogous to our iterates but with the exact  pseudo-gradient (perfect first-order feedback). In contrast to a double time-scale Tikhonov approach in \cite[Chapter 12]{FaccPang2}, $\zb(t)$ is  single time-scale, that is, it has no inner loop, and doubly regularized (by considering  $\e_t$ and $\rho_t$ as regularization parameters). Our plan now is as follows. First, we  show convergence of $\zb(t)$ to $\by(t)$ and thus, concluding convergence of $\zb(t)$ to the least-norm Nash equilibrium under a subset of conditions in Assumption~\ref{assum:parameter2}. This convergence is shown in Proposition \ref{prop:zt} under  a subset of the conditions  in Assumption~\ref{assum:parameter2}. Next, in Step 3 below we show that under the additional stepsize assumptions, our iterates $\bmu(t)$ converge to $\zb(t)$ and hence, the least-norm Nash equilibrium. 

	\begin{prop}\label{prop:zt}
		Under Assumptions~\ref{assum:convex}-\ref{assum:Lipschitz}, parts ~\ref{itm:rho_eps}),\ref{itm:omega}),\ref{itm:gam_eps}) of Assumption~\ref{assum:parameter2},  and additional condition $\sum_{t=0}^{\infty}\g_t^2 < \infty$, the sequence $\zb(t)$ converges to the least norm Nash equilibrium $\ab^*$ of the game $\Gamma(N, \{\Ab^i\}, \{J^i\})$. Moreover, if the least-norm solutions of $VI(\Mb,\Ab)$ are contained in $\mbox{int}(\Ab)$, than the condition~\ref{itm:omega}) in Assumption~\ref{assum:parameter2} can be relaxed by  the condition~2)* and the same convergence holds.
	\end{prop}
See Appendix~\ref{app:zt}    for proof. Having established the above two propositions, the main idea for the proof of Theorem \ref{th:main} is to leverage  the result from~\cite[Lemma 10]{polyak} stating almost sure convergence to 0 of a given  non-negative sequence of random numbers $a_t$, if this sequence satisfies the following condition: 
\begin{align}\label{eq:Polyak}
	&\E\{a_{t+1}|\EuScript F_t\}\le(1-\alpha_t)a_t + \phi_t,
\end{align}
where $\sum_{t=0}^{\infty}\alpha_t = \infty$,  	$\sum_{t=0}^{\infty}\phi_t < \infty$, and $\EuScript F_t$ is the $\sigma$-algebra generated by the random variables $\{a_k\}_{k\le t}$. 

First, let us demonstrate that $\|\zb(t+1) - \by(t)\|^2$ satisfies~\eqref{eq:Polyak} (its deterministic version) with $a_{t+1} = \|\zb(t+1) - \by(t)\|^2$ under the conditions of the proposition. To do so, we estimate the distance between $\zb(t+1)$ and $\by(t)$, taking into account~\eqref{eq:crmd} and the fact that $\by(t) = \Proj_{(1-\rho_t)\Ab}[\by(t) - \g_t(\Mb(\by(t))+\e_t\by(t))]$. Using non-expansion of the projection operator and conditions of the proposition, we obtain
$\|\zb(t+1) - \by(t)\|^2\le (1-\e_t\g_t)\|\zb(t) - \by(t)\|^2$. Next, to get the inequality~\eqref{eq:Polyak}, we use the relation $\|\zb(t) - \by(t)\|^2\le (1+0.5\e_t\g_t)\|\zb(t) - \by(t-1)\|^2 + \left(1+\frac{1}{0.5\e_t\g_t}\right)\| \by(t)- \by(t-1)\|^2$, as well as the estimation $\|\by(t)-\by(t-1)\|^2= O\left(\frac{|\e_t-\e_{t-1}|^2}{\e^{2}_t} + \frac{|\rho_t-\rho_{t-1}|^2}{\e^{4+2\varepsilon}_t}\boldsymbol{1}_{\ba^i\in \partial\Ab^i} \right)$ for any $\varepsilon>0$ (for the proof of the latter relations see Lemmas 8 and 9 in respectively). Finally, taking into account Assumption~\ref{assum:parameter2} 2) and 4), we obtain Inequality~\eqref{eq:Polyak} with $\alpha_t = 0.5\e_t\g_t$ and $\phi_t = \|\by(t)-\by(t-1)\|^2= O\left(\frac{|\e_t-\e_{t-1}|^2}{\e^{2}_t} + \frac{|\rho_t-\rho_{t-1}|^2}{\e^{4+2\varepsilon}_t}\boldsymbol{1}_{\ba^i\in \partial\Ab^i} \right)$ as desired:
 \begin{align*}
	\nonumber
	\|&\zb(t+1) - \by(t)\|^2\le (1-0.5\e_t\g_t)\|\zb(t) - \by(t-1)\|^2 \\
	\nonumber
	&+ O\left(\frac{|\e_t-\e_{t-1}|^2}{\e^3_t\g_t}+\frac{|\rho_t-\rho_{t-1}|^2}{\e^{5+2\varepsilon}_t\g_t}\boldsymbol{1}_{\ba^i\in \partial\Ab^i}\right).
\end{align*}



\subsubsection*{Step 3} We now analyze $\|\bmu(t)-\zb(t)\|^2$. We aim again to show that Inequality  ~\eqref{eq:Polyak} holds with $a_t = \|\bmu(t)-\zb(t)\|^2$ and
where $\EuScript F_t$, as before, is the $\sigma$-algebra generated by the random variables $\{\bmu(k),\bxi(k)\}_{k\le t}$. It then follows that $\lim_{t\to\infty}\|\bmu(t)-\zb(t)\|^2=0$ almost surely. Finally, using Proposition \ref{prop:zt} we conclude $\bmu(t)$ converges to a least-norm Nash equilibrium.

		Let us analyze each term in
		$\|\bmu(t+1)-\zb(t+1)\|^2= \sum_{i=1}^{N} \|\bmu^i(t+1)-\zb^i(t+1)\|^2$.
		From the procedures for the update of $\bmu(t)$ and $\zb(t)$ and the non-expansion property of the projection operator,
		we obtain
		\begin{align}\label{eq:nonexp1}
			\|&\bmu^i(t+1)-\zb^i(t+1)\|^2 \\
			\nonumber
			&\le \|\bmu^i(t)-\zb^i(t)\|^2 \cr
			&\quad- 2\g_t\la\Mb^i(\bmu(t))- \Mb^i(\zb(t)), \bmu^i(t)-\zb^i(t)\ra \\
			\nonumber
			&\qquad - 2\g_t\e_t\la\bmu^i(t)-\zb^i(t), \bmu^i(t)-\zb^i(t)\ra\cr
			&\qquad-2\g_t\la\Qb^i(t)+\Rb^i(t), \bmu^i(t)-\zb^i(t)\ra+ \g^2_t\|\Gb^i(t)\|^2,
		\end{align}
		where, for ease of notation, we have defined
		$\Gb^i(t) = \e_t(\bmu^i(t)-\zb^i(t)) + \Mb^i(\bmu(t)) - \Mb^i(\zb(t))
			+\Qb^i(t)+\Rb^i(t)$.
			It follows that the terms in $\|\Gb^i(t)\|^2$ are given as
		\begin{align*}
			&\|\Gb^i(t)\|^2 = {\epsilon^2(t)\|\bmu^i(t)-\zb^i(t)\|^2} \cr
			&\qquad+ \|\Mb^i(\bmu(t)) - \Mb^i(\zb(t))\|^2+ \|\Qb^i(t)\|^2+\|\Rb^i(t)\|^2\cr
			&\qquad+ 2\la\Qb^i(t),\Rb^i(t)\ra \cr
			&\qquad+ 2\e_t\la\Mb^i(\bmu(t)) - \Mb^i(\zb(t)),\bmu^i(t)-\zb^i(t)\ra \cr
			&\qquad+2\e_t\la\bmu^i(t)-\zb^i(t),\Qb^i(t)+\Rb^i(t)\ra \cr
			&\qquad+ 2\la\Mb^i(\bmu(t)) - \Mb^i(\zb(t)),\Qb^i(t)+\Rb^i(t)\ra.
		\end{align*}
		Thus, accounting for the above,  for Proposition~\ref{prop:sample_grad}, which implies $\E\{\Rb^i(t)|\EuScript F_t\} = 0$ for any $t$, and for the Cauchy-Schwarz inequality, we get the following from Inequality \eqref{eq:nonexp1}:
		\begin{align*}
			&\E\{\|\bmu^i(t+1)-\zb^i(t+1)\|^2|\EuScript F_t\} \cr
			&\le 
                (1-2\g_t\e_t)\|\bmu^i(t)-\zb^i(t)\|^2 \cr
			&\quad- 2\g_t\la\Mb^i(\bmu(t))- \Mb^i(\zb(t)), \bmu^i(t)-\zb^i(t)\ra \cr
			&\quad + 2\g_t\E\{\|\Qb^i(t)\|\|\bmu^i(t)-\zb^i(t)\||\EuScript F_t\} \cr
			&\quad+ \g^2_t\epsilon^2_t\|\bmu^i(t)-\zb^i(t)\|^2 + \g^2_t\cr
			&\quad+ \g^2_t\E\{\|\Qb^i(t)\|^2|\EuScript F_t\}+\g^2_t\E\{\|\Rb^i(t)\|^2|\EuScript F_t\}\cr
			&\quad+2\g^2_t\E\{\|\Qb^i(t)\|\|\Rb^i(t)\||\EuScript F_t\}\cr
			&\quad +2\g^2_t\e_t\|\bmu^i(t)-\zb^i(t)\|\E\{\|\Qb^i(t)\|\|\Rb^i(t)\||\EuScript F_t\} \cr
			&\quad+ 2\g_t^2\|\Mb^i(\bmu(t)) - \Mb^i(\zb(t))\|\E\{\|\Qb^i(t)\||\EuScript F_t\},
		\end{align*}
		where in the last inequality we used $\|\Mb^i(\bmu(t)) - \Mb^i(\zb(t))\|^2 = O(\|\bmu(t) - \zb(t)\|^2)$ implied by Assumptions~\ref{itm:lip} and~\ref{assum:Lipschitz}  and Remark~\ref{rem:LipschJ}.
		Next, by taking into account compactness of the set $\Ab$ and the properties of the terms~$\Qb^i$ and $\Rb^i$ (see Proposition \ref{prop:sample_grad}), we conclude that
		\begin{align*}
			&\E\{\|\bmu^i(t+1)-\zb^i(t+1)\|^2|\EuScript F_t\} \cr
			& \le(1-2\g_t\e_t)\|\bmu^i(t)-\zb^i(t)\|^2 \cr
			&\quad- 2\g_t\la\Mb^i(\bmu(t))- \Mb^i(\zb(t)), \bmu^i(t)-\zb^i(t)\ra \cr
			&\quad+O(\g_t\s_t + \g_t^2\e_t^2 + \g_t^2 + \frac{\g_t^2}{\s_t^2}).
		\end{align*}
		As the pseudo-gradient is monotone, $\la\Mb(\bmu(t))- \Mb(\zb(t)), \bmu(t)-\zb(t)\ra\ge 0$. Hence, as $\e_t\to 0$, from the inequalities above for $i=1,\ldots,N$ we conclude 
		\begin{align}\label{eq:final_main}
			\E\{&\|\bmu(t+1)-\zb(t+1)\|^2|\EuScript F_t\}\cr
			\le&(1-2\e_t\g_t)\|\bmu(t) - \zb(t)\|^2 + O(\g_t\s_t +  \g_t^2 + \frac{\g_t^2}{\s_t^2}).
		\end{align}
		From Assumption~\ref{assum:parameter2} we get
		$\sum_{t=0}^\infty(\g_t\s_t +  \g_t^2 + \frac{\g_t^2}{\s_t^2}) < \infty$, $\sum_{t=0}^\infty \e_t\g_t=\infty$.  Thus, \eqref{eq:Polyak} holds and consequently, $\|\bmu(t)-\zb(t)\|^2$  converges almost surely to $0$ as $t\to\infty$. Next, by taking into account Proposition~\ref{prop:zt}, we obtain that
		$\Pr\{\lim_{t\to\infty}\bmu(t)=\ab^*\}=1,$
		where $\ab^*$ is the least-norm Nash equilibrium in the game $\Gamma(N, \{\Ab^i\}, \{J^i\})$.
		Finally, since $\sigma_t\to 0$, we can apply Portmanteau lemma \cite{portlem} to conclude that $\ab(t)$ convergence in probability to $\ab^*$ as desired.
\end{proof}

\begin{rem}
\label{rem:reg_rate}
The proof of Theorem~\ref{th:main} implies the rate with which the sequence $\bmu(t)$ converges to $\zb(t)$, by applying for example,  Chung's lemma~\cite{Chung} to Inequality~\eqref{eq:final_main}. Similarly, we can  establish the convergence rate of $\zb(t)$ to $\by(t)$  defined in~\eqref{eq:yt}. However, the rate at which $\by(t)$ converges to the least norm solution of the variational inequality $VI(\Ab, \Mb)$ can be obtained under additional assumptions  such as linear mappings with ``source representable" solutions, see \cite[Chapter 5.8]{bakushinsky2012ill}.  {Note that the setting of purely monotone games under consideration corresponds to convex (non-strictly) optimization. In the latter,  in general we can have convergence rate in terms of function values but not in terms of the algorithm iterates}. 
\end{rem}

 {Next, we consider the performance of the algorithm  in terms of the regret benchmark. We characterize the  regret rate in terms of the algorithm parameters. This characterization can guide a player in setting the parameters of her algorithm.} 

\section{No-regret Property of the Algorithm}\label{sec:regret}
We first interpret the problem of learning under payoff information from the perspective of a single player as an instance of online convex optimization.  Next,  we prove the no-regret property of  Algorithm \ref{alg:payoff} and derive its regret rate.

In an online convex optimization (or online convex programming), at each time step $t$,  each player $i$ selects an action $\ba^i(t) \in \Ab^i$. After the selection, she receives the cost function $ J^i_t=J^i_t(\ba^i(t)) = J^i(\ba^i(t), \ba^{-i}(t))$. Thus, the player experiences a sequence of cost functions which are a priori unknown to her due to their dependence on  actions of other players. Let us denote this sequence by   $\{J^i_1, \ldots, J^i_t, \ldots \}$, where $J^i_t(\cdot) = J^i(\cdot, \ba^{-i}(t))$, $J_t^i : \R^{d} \to \R$. Note that each $J^i_t$ is  convex  on $\R^d$ and smooth due to Assumptions ~\ref{assum:convex} and ~\ref{assum:Lipschitz}, respectively.} Furthermore, consistent with the payoff information setup of the game, each player at each time $t$ receives only the cost of the played action, $J^i_t(\ab^i(t))$, as  feedback, rather than the full cost function $J^i_t(.)$ or its gradient. {This setting is referred to as bandit online convex optimization \cite{flaxman2005online, saha2011improved}}.

In an online convex optimization problem, efficiency of  player $i$'s algorithm is measured with respect to its regret. 			
\begin{definition}\label{def:regret}
Given an algorithm updating $\{\ba^i(t)\}$ with respect to an online convex programming problem $(\Ab^i,\{J^i_1, J^i_2,\ldots\})$, if $\{\ba^i(1), \ba^i(2),\ldots\}$ are the vectors selected by this algorithm,  the \emph{regret} of the algorithm until time $T$ with respect to any reference point $\ab^i\in\Ab^i$ is
\begin{align}
\label{eq:regret_eq}
R^i(T,\ab^i) = \sum_{t=1}^T J^i_t(\ba^i(t)) -  \sum_{t=1}^T  J^i_t (\ba^i ).
\end{align}
\end{definition}

The goal of an online optimization algorithm is to guarantee  \emph{no-regret}. That is,
		\begin{align*}
			\varlimsup_{T\to\infty}\frac{R^i(T,\ab^i)}{T}=\lim_{T\to\infty} \sup_T\frac{R^i(T,\ab^i)}{T}\le 0.
		\end{align*}

	
	\begin{theorem}\label{th:no-regret_regular}
		For a given player $i$, choose the parameters in  Algorithm \ref{alg:algorithm1} as
		$\gamma_t = O\left(\frac{1}{t^g}\right)$, $\sigma_t = O\left(\frac{1}{t^s}\right)$, $\epsilon_t= O\left(\frac{1}{t^e}\right)$, $\rho_t = O\left(\frac{1}{t^r}\right)$,  with $ 0< g,s,r<1$, and 
\begin{enumerate} 
\item $r<s$. 
\item $g-2s>0$.
\end{enumerate}
Let other players choose their actions arbitrarily. Under Assumptions~\ref{assum:convex}-\ref{assum:Lipschitz} and given an arbitrary $\bmu^i(0)$ Algorithm \ref{alg:algorithm1} is no-regret for player $i$ with
		\begin{align}
			\label{eq:reg_rate}
			&\E \frac{R^i(T,\ab^i)}{T} \cr
			&=  O (\frac{1}{T^{1-g }} + \frac{1}{T^{g-2s}} + \frac{1}{T^{2s}}+ \frac{1}{T^{r}}\boldsymbol{1}_{\ab^i\in \partial \Ab^i} + \frac{1}{T^{e}}),
		\end{align}
  where the expectation is taken in respect to the sequence of stochastic vectors $\{\ba^i(t)\}$.
	\end{theorem}

\begin{proof}
		First, we connect the regret $R^i(T,\ab^i)$ defined in \eqref{eq:regret_eq} with a regret based on the mean vectors $\bmu^i(t)$, $t=1,2,\ldots$. Let us introduce the following notation:
		\begin{align}\label{eq:Jtmu}
			J^i_{t,\bmu}(\cdot) = J^i(\cdot,\bmu^{-i}(t)).
		\end{align}
		Given  $\ab^i\in \Ab^i$ define the $\bmu$-regret $R_{\bmu}^i(T,\ab^i)$ as:
		$$R_{\bmu}^i(T,\ab^i) = \sum_{t=1}^T [J^i_{t,\bmu}(\bmu^i(t)) - J^i_{t,\bmu}(\ba^i)].$$
		The lemma below connects $R_{\bmu}^i(T,\ab^i)$ and $R^i(T,\ab^i)$.
		\begin{lem}\label{lem:muRegret}
			Choose the parameters  $\sigma_t = \frac{1}{t^s}$ and $\rho_t = \frac{1}{t^r}$  with $0<s,r$, $r<s$.
			Then, under Assumption~\ref{assum:Lipschitz}~\ref{itm:lip}),
			\[\E R^i(T,\ab^i) = \E R_{\bmu}^i(T,\ab^i) + O(\sigma_t^2),\]
   where the expectation is taken in respect to the sequence  $\{\ba^i(t)\}_t$ (left hand side) and $\{\bmu^i(t)\}_t$ (right hand side).
		\end{lem}
		Please see Appendix~\ref{app:muRegret} for the proof.
		
		Now we focus on the estimation of the expected value of the $\bmu$-regret, namely $\E R_{\bmu}^i(T,\ab^i)$, for the procedure defined by Algorithm~\ref{alg:algorithm1}. Due to definition of $\Qb^i$ in Equation~\eqref{eq:Qi0}, the procedure~\eqref{eq:rppg_i} for each player $i$ can be rewritten as follows:
		\begin{align}
			\label{eq:pbavmuOL0}
			\bmu^i&(t+1) =\Proj_{(1-\rho_t)\Ab^i}[\bmu^i(t) \cr
			&-\gamma_t( \tilde{\Mb}^i(\bmu(t))+ \Rb^i(t) + \Pb^i(t)+\e_t\bmu^i(t))],
		\end{align}
		where $\tilde{\Mb}^i$, $\Rb^i$ were defined in \eqref{eq:smoothed_main} and \eqref{eq:Ri_main}, whereas 
		$\Pb^i(t) = J^i(\ab(t))-J^i(\bxi(t)))\frac{\bxi^i(t) -\bmu^i(t)}{\sigma^2_t}$
		is term 2 in the definition of $\Qb^i$ (see Equation ~\eqref{eq:Qi0}).
		Next, let the function $\tilde{J}_{t,\bmu}^i$ be defined as follows:
		\begin{align}\label{eq:tildeJ_t}
			\tilde{J}_{t,\bmu}^i(\bmu^i) = \int_{\R^{Nd}}J^i(\bx)p^i(\bx^i;\bmu^i,\sigma_t)\prod_{j\ne i}p^j(\bx^j;\bmu^j(t),\sigma_t)d\bx,
		\end{align}
		where $p^k(\bx^k;\bmu^k,\sigma_t)$, $k=1,\ldots,N$, is the density function defined in Equation~\eqref{eq:density}.
		Then \eqref{eq:pbavmuOL0} can be written as
		\begin{align}
			\label{eq:pbavmuOL}
			\bmu^i&(t+1) =\Proj_{(1-\rho_t)\Ab^i}[\bmu^i(t) \cr&-\gamma_t(\nabla\tilde{J}_{t,\bmu}^i(\bmu^i(t)) + \Rb^i(t)+ \Pb^i(t) +\e_t\bmu^i(t))],
		\end{align}
		where we used Equations~\eqref{eq:smoothed_main} and~\eqref{eq:tildeJ_t} to get
		\begin{align*}
			&\tilde{\Mb}^i(\bmu(t)) =\frac{\partial \tilde{J}^i(\bmu(t))}{\partial \bmu^i} \cr
			&= \frac{\int_{\R^{Nd}}J^i(\bx)p^i(\bx^i;\bmu^i(t),\sigma_t)\prod_{j\ne i}p^j(\bx^j;\bmu^j(t),\sigma_t)d\bx}{\partial \bmu^i}\cr
			&\qquad\qquad\qquad = \nabla\tilde{J}_{t,\bmu}^i(\bmu^i(t)).
		\end{align*}
		Next, we provide a result connecting the function $J_{t,\bmu}^i(\bmu^i)$ with its mixed strategy version $\tilde{J}_{t,\bmu}^i(\bmu^i)$.
		\begin{lem}
			\label{lem:mixedstr_cost}
			Under Assumption~\ref{assum:convex} and part ~\ref{itm:lip} of Assumption~\ref{assum:Lipschitz}, the functions $\tilde{J}_{t,\bmu}^i:\R^d \rightarrow \R$ satisfies the following:
			\begin{enumerate}
				\item \label{itm:convex_smooth} $\tilde{J}_{t,\bmu}^i$, $t=1,\ldots$, are convex on $\R^d$ and their gradients $\nabla \tilde{J}_{t,\bmu}^i$ are bounded over $\Ab$.
				\item \label{itm:bound_smooth}
					$|J_{t,\bmu}^i(\bmu^i)-\tilde{J}_{t,\bmu}^i(\bmu^i)|=O(\sigma_t^2), \; \forall \bmu^i \in \Ab^i.$
			\end{enumerate}
		\end{lem}
	Note that Lemma~\ref{lem:mixedstr_cost} has been derived in \cite{NesterovSpokoiny} for the Gaussian smoothing and in~\cite{PerchetBach} in the case of the smoothing technique based on the uniform distribution over the unit sphere. We provide a proof in Appendix~\ref{app:mixed_cost} for completeness. 
	We now proceed further with the proof of Theorem \ref{th:no-regret_regular}. 
 
 Let us consider any $\ba^i\in \partial\Ab^i$. Then,
	\begin{align}\label{eqR:withrt0}
	&\tilde{J}_{t,\bmu}^i(\bmu^i(t))- \tilde{J}_{t,\bmu}^i(\ba^i) = \tilde{J}_{t,\bmu}^i(\bmu^i(t)) \\- \; &\tilde{J}_{t,\bmu}^i(\Proj_{(1-\rho_t)\Ab^i}\ba^i) 
		 + \tilde{J}_{t,\bmu}^i(\Proj_{(1-\rho_t)\Ab^i}\ba^i) - \tilde{J}_{t,\bmu}^i(\ba^i).
    \nonumber
	\end{align}
	Due to the compactness of $\Ab^i$ and smoothness of $\tilde{J}_{t,\bmu}^i$ from Lemma~\ref{lem:mixedstr_cost} part ~\ref{itm:convex_smooth}, $\tilde{J}_{t,\bmu}^i$ is Lipschitz continuous over $\Ab$. Thus,
	\[\tilde{J}_{t,\bmu}^i(\Proj_{(1-\rho_t)\Ab^i}\ba^i) - \tilde{J}_{t,\bmu}^i(\ba^i)\le \tilde L^i\|\Proj_{(1-\rho_t)\Ab^i}\ba^i - \ba^i\|,\]
	for some positive constant $\tilde L^i$.
	Taking into account that for any $\ba^i\in \Ab^i$ 
	$
		\|\Proj_{(1-\rho_{t_1})\Ab^i}\ba^i - \Proj_{(1-\rho_{t_2})\Ab^i}\ba^i\| = O(|\rho_{t_1} - \rho_{t_2}|),
	$
	(see Lemma~7 in Appendix~D  for an analogous proof) we get from Equation \eqref{eqR:withrt0}
	$	
	\tilde{J}_{t,\bmu}^i(\bmu^i(t))- \tilde{J}_{t,\bmu}^i(\ba^i)\le  \tilde{J}_{t,\bmu}^i(\bmu^i(t)) - \tilde{J}_{t,\bmu}^i(\Proj_{(1-\rho_t)\Ab^i}\ba^i) + O(\rho_t).
	$
	On the other hand, if $\ba^i\notin \partial\Ab^i$, let us choose $\rho_t$ such that for all $t$ we have
	$\ba^i\in (1-\rho_t) \Ab^i$ and, thus, $\ba^i=\Proj_{(1-\rho_t)\Ab^i}\ba^i$. Hence,
	\begin{align}\label{eqR:eq}
	&\tilde{J}_{t,\bmu}^i(\bmu^i(t))- \tilde{J}_{t,\bmu}^i(\ba^i)\cr&\le \tilde{J}_{t,\bmu}^i(\bmu^i(t)) - \tilde{J}_{t,\bmu}^i(\Proj_{(1-\rho_t)\Ab^i}\ba^i) + O(\rho_t)\boldsymbol{1}_{\ba^i\in \partial\Ab^i}.
	\end{align}
	Moreover, due to Lemma~\ref{lem:mixedstr_cost}, the function $\tilde{J}_{t,\bmu}^i$ defined by Equation~\eqref{eq:tildeJ_t} is convex. Thus, for any $\ab^i\in \Ab^i$,
	\begin{align}\label{eqR:eq1}
		\tilde{J}_{t,\bmu}^i&(\bmu^i(t)) - \tilde{J}_{t,\bmu}^i(\Proj_{(1-\rho_t)\Ab^i}\ba^i)\cr&\le \la\nabla \tilde{J}_{t,\bmu}^i(\bmu^i(t)), \bmu^i(t)-\Proj_{(1-\rho_t)\Ab^i}\ba^i\ra.
	\end{align}
	Finally, due to the fact that $|J_{t,\bmu}^i(\ba^i)-\tilde{J}_{t,\bmu}^i(\ba^i)|=O(\sigma_t^2)$, $\forall \ba^i \in \Ab^i$, (see Lemma~\ref{lem:mixedstr_cost} part~\ref{itm:bound_smooth}), we obtain
	\begin{align}
		\label{eqR:reg_connect}
		R_{\bmu}&(T,\ab^i) = \sum_{t=1}^T  J_{t,\bmu}^i(\bmu^i(t)) - \sum_{t=1}^T J_{t,\bmu}^i(\ab^i) \cr&= \sum_{t=1}^T(\tilde{J}_{t,\bmu}^i(\bmu^i(t))-\tilde{J}_{t,\bmu}^i(\ab^i)) + \sum_{t=1}^T O(\sigma_t^2).
	\end{align}
	Thus, given any reference point $\ba^i \in \Ab^i$, in view of Inequalities~\eqref{eqR:eq}, \eqref{eqR:eq1}, and \eqref{eqR:reg_connect}, we can write
	\begin{align}\label{eqR:final}
		&R_{\bmu}(T,\ab^i) \le \sum_{t=1}^T\la\nabla \tilde{J}_{t,\bmu}^i(\bmu^i(t)), \bmu^i(t)-\Proj_{(1-\rho_t)\Ab^i}\ba^i\ra\cr
  &\qquad+ \sum_{t=1}^TO(\rho_t)\boldsymbol{1}_{\ba^i\in \partial\Ab^i}+\sigma^2_t).
	\end{align}
We  focus further on the term $\la\nabla \tilde{J}_{t,\bmu}^i(\bmu^i(t)), \bmu^i(t)-\Proj_{(1-\rho_t)\Ab^i}\ba^i\ra$.
In the following analysis, for the sake of presentation clarity, we  use the following notations:  $\Proj_{(1-\rho_t)\Ab^i}\ba^i = \ba_{P,t}^i$.
	Due to non-expansion of the projection operator and Equation~\eqref{eq:pbavmuOL}, we get for any $\ba^i\in \Ab^i$
 \begin{align}\label{eq:regret_0}
		&\|\bmu^i(t+1) - \ba_{P,t}^i\|^2\le\|\bmu^i(t) - \ba_{P,t}^i\|^2 \cr
				&+ \gamma^2_t(\|\nabla\tilde{J}_{t,\bmu}^i (\bmu^i(t)) \|^2 +\|\Rb^i(t)\|^2+\|\Pb^i(t)\|^2+\epsilon^2_t\|\bmu^i(t)\|^2)\cr
			&-2\gamma_t\la\nabla\tilde{J}_{t,\bmu}^i (\bmu^i(t)),\bmu^i(t)-\ba_{P,t}^i\ra\cr
		&-2\gamma_t\la\Rb^i(t)+ \Pb^i(t)+\epsilon_t\bmu^i(t),\bmu^i(t)-\ba_{P,t}^i\ra\cr
			&+ 2\gamma^2_t[\la\nabla\tilde{J}_{t,\bmu}^i (\bmu^i(t)),\Pb^i(t)\ra\cr
		&+ \la\nabla\tilde{J}_{t,\bmu}^i (\bmu^i(t)),\Rb^i(t)\ra + \la\nabla\tilde{J}_{t,\bmu}^i (\bmu^i(t)),\epsilon_t\bmu^i(t)\ra  \cr
			&+ \la\Rb^i(t),\epsilon_t\bmu^i(t)\ra+\la\Pb^i(t),\Rb^i(t)\ra+\la\Pb^i(t),\epsilon_t\bmu^i(t)\ra].
	\end{align}
	Thus, taking into account the properties of the term  $\Rb^i(t)$ from Proposition~\ref{prop:sample_grad}, boundedness of $\|\nabla\tilde{J}_{t,\bmu}^i (\bmu^i(t)) \|$ from Lemma~\ref{lem:mixedstr_cost} part~\ref{itm:convex_smooth}), and~\eqref{eq:regret_0}, we obtain that for $r<s$,
	\begin{align*}
		&\E\{\|\bmu^i(t+1) - \ba_{P,t}^i\|^2|\EuScript F_t\}\le\|\bmu^i(t) - \ba_{P,t}^i\|^2\cr
		&\quad-2\gamma_t\la\nabla\tilde{J}_{t,\bmu}^i (\bmu^i(t)),\bmu^i(t)-\ba_{P,t}^i\ra+ O\left(\frac{\gamma^2_t}{\sigma_t^2}\right) + O(\gamma_t\e_t).
	\end{align*}
	Hence, by taking the full expectation of the both sides and rearranging the terms, we get
	\begin{align*}
		&\E\la\nabla\tilde{J}_{t,\bmu}^i (\bmu^i(t)),\bmu^i(t)-\ba_{P,t}^i\ra\le  O\left(\frac{\gamma_t}{\s_t^2}+\epsilon_t\right)\cr& \frac{1}{2\gamma_t}[\E\|\bmu^i(t) - \ba_{P,t}^i\|^2-\E\|\bmu^i(t+1) - \ba_{P,t}^i\|^2]
	\end{align*}
	By summing up the left and the right hand side over $t$ and taking into account that $\|\bmu^i(t+1) - \ba_{P,t}^i\|^2\le 2M$ for some $M>0$, we obtain
	\begin{align*}
		&\sum_{t=1}^T \E\la\nabla\tilde{J}_{t,\bmu}^i (\bmu^i(t)),\bmu^i(t)-\ba_{P,t}^i\ra\le \frac{1}{2\gamma_1}\E\|\bmu^i(1) - \ba_{P,t}^i\|^2 \cr
		&+ \sum_{t=2}^T \E\|\bmu^i(t+1) - \ba_{P,t}^i\|^2\left(\frac{1}{2\gamma_t } - \frac{1}{2\gamma_{t-1} }\right) + \sum_{t=1}^TO\left(\frac{\gamma_t}{\s_t^2}\right)\cr
		&\qquad+ \sum_{t=1}^TO(\epsilon_t) \le \frac{M}{\gamma_T }  + \sum_{t=1}^TO\left(\frac{\gamma_t}{\s_t^2}\right)+ \sum_{t=1}^TO(\epsilon_t).
	\end{align*}
	Thus, due to~\eqref{eqR:final},
	\begin{align*}
		\E R_{\bmu}(T,\ab^i)& \le \frac{M}{\gamma_T }  + \sum_{t=1}^TO\left(\frac{\gamma_t}{\s_t^2}\right)+\sum_{t=1}^TO(\epsilon_t) + \sum_{t=1}^T O(\sigma^2_t)\cr &\qquad+\sum_{t=1}^T O(\rho_t)\boldsymbol{1}_{\ba^i\in \partial\Ab^i}.
	\end{align*}
	As
	$\sum_{t=1}^T{\frac{\gamma_t}{\s_t^2}}= O(T^{1-g+2s})$, $\sum_{t=1}^T{\rho_t} = O(T^{1-r})$, $\sum_{t=1}^T{\epsilon_t} = O(T^{1-e})$,  $\sum_{t=1}^T\sigma^2_t = O(T^{1-2s})$,
	we conclude that
	$
		\E \frac{R_{\bmu}(T,\ab^i)}{T} = O \left(\frac{1}{T^{1-g }} + \frac{1}{T^{g-2s}} + \frac{1}{T^{2s}}+ \frac{1}{T^{r}}\boldsymbol{1}_{\ba^i\in \partial\Ab^i} + \frac{1}{T^{e}}\right).
	$
	Finally, applying  Lemma~\ref{lem:muRegret}, we obtain the result.
\end{proof}

Note that any parameters satisfying Assumption~\ref{assum:parameter2} and chosen based on the p-series conditions in Lemma~\ref{lem:step_constants} ensure the parameter conditions in the theorem above are satisfied. Hence, the algorithm converging to a least-norm Nash equilibrium is also no-regret. Furthermore, $\epsilon_t$ can be arbitrarily small and the algorithm will still be no-regret. However, for convergence to a Nash equilibrium, $\epsilon_t$ should be at most $O(1/t)$.   {The theorem above  extends the class of games for which a no-regret payoff-based algorithm has provable convergence to a Nash equilibrium, from potential  \cite{heliou2017learning} and strongly monotone \cite{bravo2018bandit}  to merely monotone games}. 

Next, what choices of parameters  would optimize the regret rate? Below, we address this question. Below, we address this question in two  cases: 1) a single player implementing the algorithm; 2) all players implementing the algorithm.  
\begin{cor}
\label{cor:regret_rate}
Let  $\varepsilon > 0$ denote a small constant.

\begin{enumerate}
\item If player $i$ implements Algorithm \ref{alg:algorithm1} with an arbitrary $\bmu^i(0)$, then regardless of other players choice of algorithm, her optimal regret rate with respect to $\ba^i\in \mbox{int}(\Ab^i)$ is  $O\left(\frac{1}{T^{1/3}}\right)$ and with respect to $\ba^i\in \partial\Ab^i$ is $O\left(\frac{1}{T^{1/4-\varepsilon}}\right)$;
\item If all players implement Algorithm \ref{alg:algorithm1}, and the parameters satisfy \ref{assum:parameter2} for convergence condition to a Nash equilibrium, the optimal regret rate with respect to $\ba^i\in \mbox{int}(\Ab^i)$ is  $O\left(\frac{1}{T^{1/4-\varepsilon}}\right)$ and  {with respect to $\ba^i\in \partial\Ab^i$ is $O\left(\frac{1}{T^{1/6-\varepsilon}}\right)$}.
 
\label{cor:rate2}
\end{enumerate}

\end{cor}
\begin{proof}
The optimal rates are obtained by maximizing the minimum exponent of $T$, in the four terms of \eqref{eq:reg_rate} under the condition of  Theorem \ref{th:no-regret_regular} in case 1) and Theorem \ref{th:main} and Lemma \ref{lem:step_constants} in case 2).  {This maximization can be done by solving a linear program as the objectives and constraints on the parameters $(g,s,e,r)$ are linear.} {The setting of the  parameters achieving the above optimal rates are provided  in Table \ref{table:rates}. In the table, $\varepsilon > 0$ is an arbitrary small constant and  $0 < g, s, r, e < 1$. }


\begin{table}
\caption{Setting of $(g,s,r,e)$ for optimal regret rate. }
\label{table:rates}
\centering
\begin{tabular}{ |c|l|l| } 
\hline
& regret  & Nash 
\\
\hline
interior & $(\frac{2}{3},\frac{1}{6}, < \frac{1}{6}, \geq \frac{1}{3})$ & $(\frac{3}{4} + \varepsilon,\frac{1}{4}, \frac{1}{4}-\frac{\varepsilon}{2}, \frac{1}{4}-\varepsilon)$ \\ 
\hline
boundary & $(\frac{3}{4},\frac{1}{4}, \frac{1}{4}-{\varepsilon}, \geq \frac{1}{4})$ & $(\frac{5}{6} +\varepsilon, \frac{1}{3}, \frac{1}{3}-{\varepsilon}, \frac{1}{6}-\varepsilon)$ 
\\
\hline
\end{tabular}

\end{table}


\end{proof}
		
\begin{rem}
\label{rem:reg_rate}
The rate $O\left(\frac{1}{T^{1/4}}\right)$ is consistent with the regret rate derived in \cite{flaxman2005online} for the online learning of the class of smooth convex  functions. This rate was improved to  $O\left(\frac{1}{T^{1/3}}\right)$ in \cite{saha2011improved}. We have been able to achieve this improved regret rate only for an interior reference point. 
Furthermore,  requiring convergence to a Nash equilibrium worsens the regret rate of the algorithm. We believe the above regret rates are suboptimal. Future research can explore other classes of payoff-based learning algorithms for monotone games to potentially  improve these rates.
			
	\end{rem}
	


\section{Simulations}\label{sec:sim}

We illustrate the performance of the algorithm by addressing Examples~\ref{example:penny} and~\ref{example:nonunique} from Section~\ref{sec:example}.  {We further illustrate the applicability of the algorithm to higher dimensional problems.}
The parameters of our algorithm were set to satisfy convergence condition in Theorem~\ref{th:main} as per Lemma~\ref{lem:step_constants}, while being  optimal for minimizing regret as per Corollary~\ref{cor:regret_rate} ($\varepsilon$ was set to  $0.04$, see Table \ref{table:rates}). In particular, they were set to    $(g,s,r,e)=(0.79, 0.25, 0.23, 0.21 )$ and $(g, s, r, e)=(0.87, 0.33, 0.29, 0.13)$, for Nash equilibria in the interior and boundary of the set $\Ab$, respectively. 


\begin{figure}[!t]
\begin{subfigure}
		\centering
		\includegraphics[width = \linewidth]{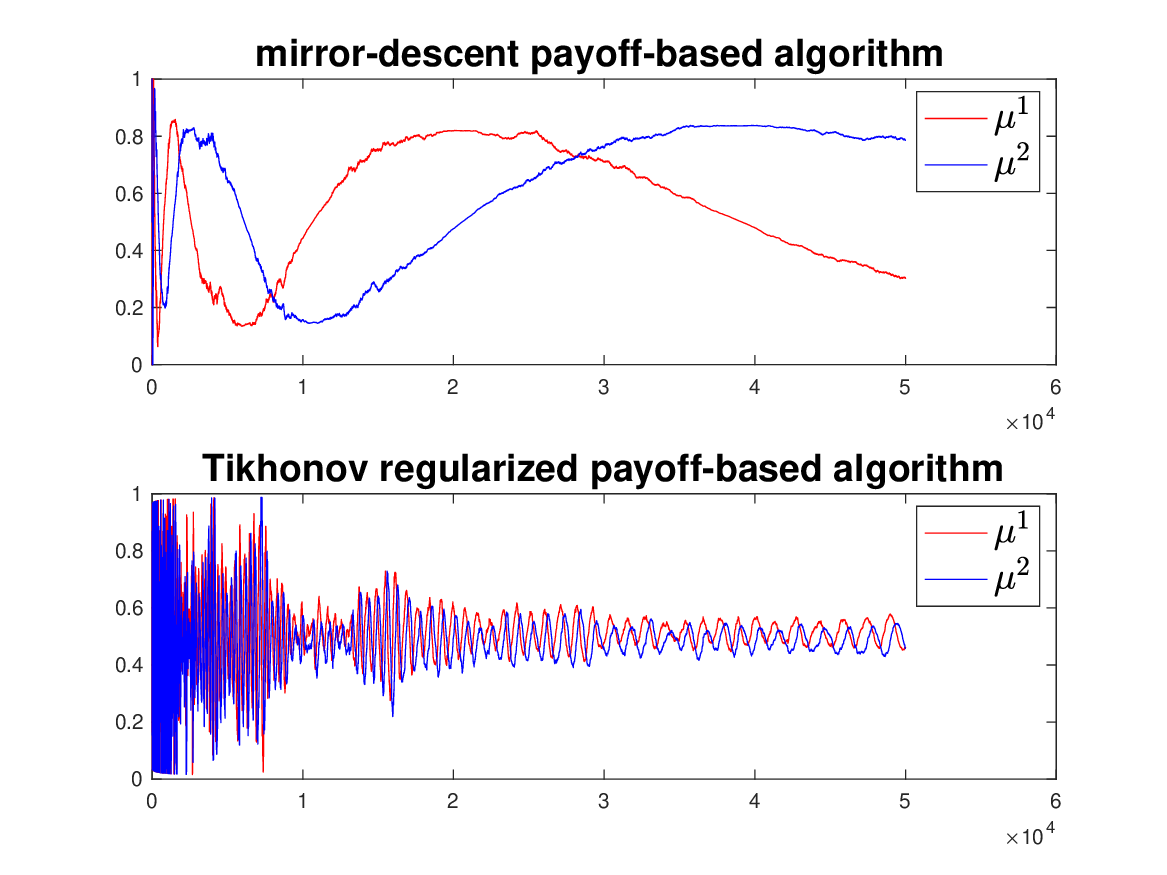}
		\caption{The iterates $\bmu(t)$ in Example \ref{example:penny}a) with $\Ab^1=\Ab^2=[0,1]^2$.}
		\label{fig:penny}
	\end{subfigure}%
\end{figure}

    \begin{figure}[!t]
	\begin{subfigure}
		\centering
		\includegraphics[width=\linewidth]{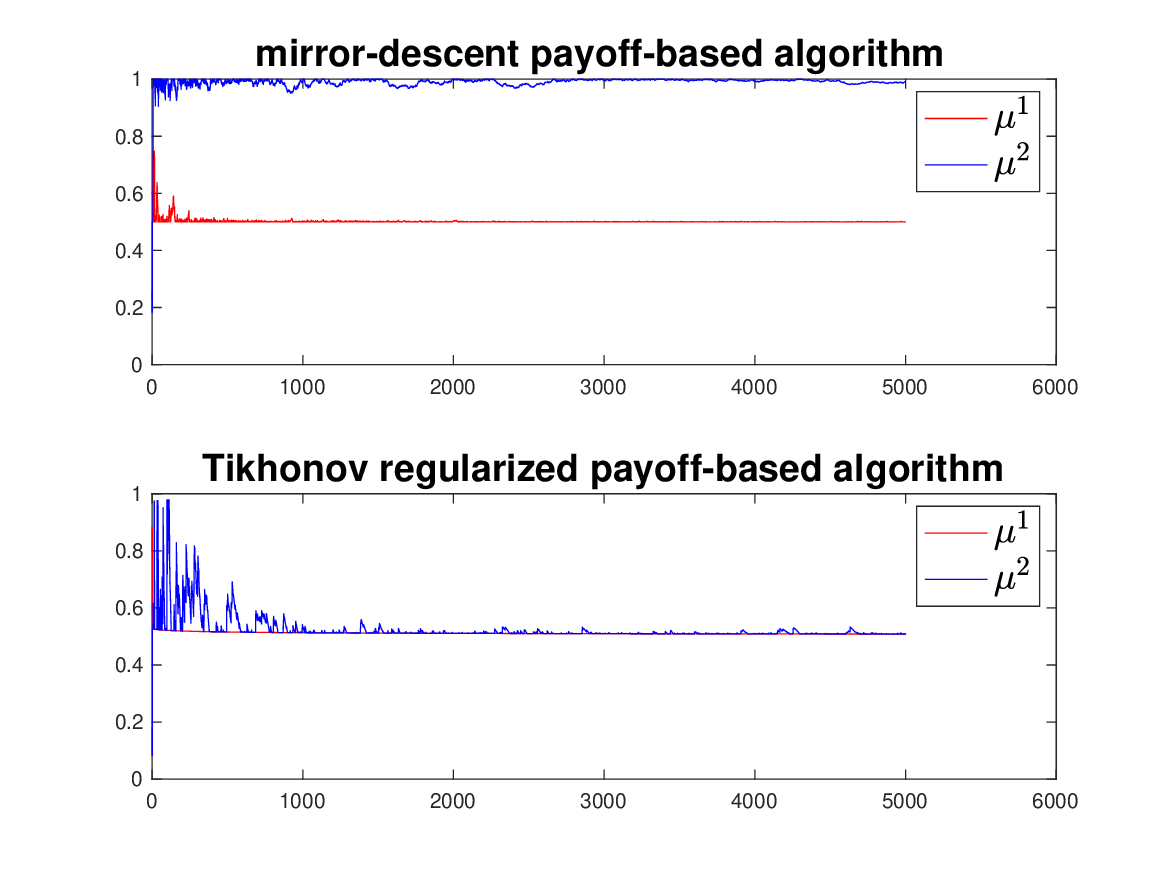}
		\caption{The iterates $\bmu(t)$ in Example \ref{example:penny}b) with $\Ab^1=\Ab^2=[1/2,1]$.}
		\label{fig:restricted}
		\end{subfigure}
\end{figure}

\subsection{Example 1} We implement the proposed algorithm for Examples 1a and 1b. Figure~\ref{fig:penny} and ~\ref{fig:restricted} illustrate the iterates of the mirror-descent algorithm  from \cite{bravo2018bandit} (top panels) and the proposed algorithm (bottom panels). The mirror-descent  algorithm \cite{bravo2018bandit} is closest to our setting since it is payoff-based, ensures feasibility of actions and is no-regret. We set its parameters as  $p = 1, q = 1/3$ (in accordance to Theorem 5.2 in \cite{bravo2018bandit}).  In Example 1a, Algorithm \ref{alg:algorithm1}'s iterate converge to the unique Nash equilibrium while in Example 1b, the iterates converge to the least-norm Nash equilibrium. Note that in a merely monotone setting, unlike the strictly monotone case, the mirror-descent algorithm does not have guaranteed convergence. In Example 1a the iterates do not converge while in 1b they  appear to converge to a non least-norm Nash equilibrium\footnote{ {In this particular example, convergence could be intuitively explained by considering mirror descent under  first-order feedback.  In this case, as $a^i \geq 1/2$, player 1's gradient is positive, whereas player 2's gradient is negative. Thus, the iterates of mirror descent approach the equilibrium $(1,1/2)$ and it can be verified that they are the fixed point of the mirror descent. }}

\begin{figure}[!t]
\centering
\begin{subfigure}
\centering
\includegraphics[width=\linewidth]{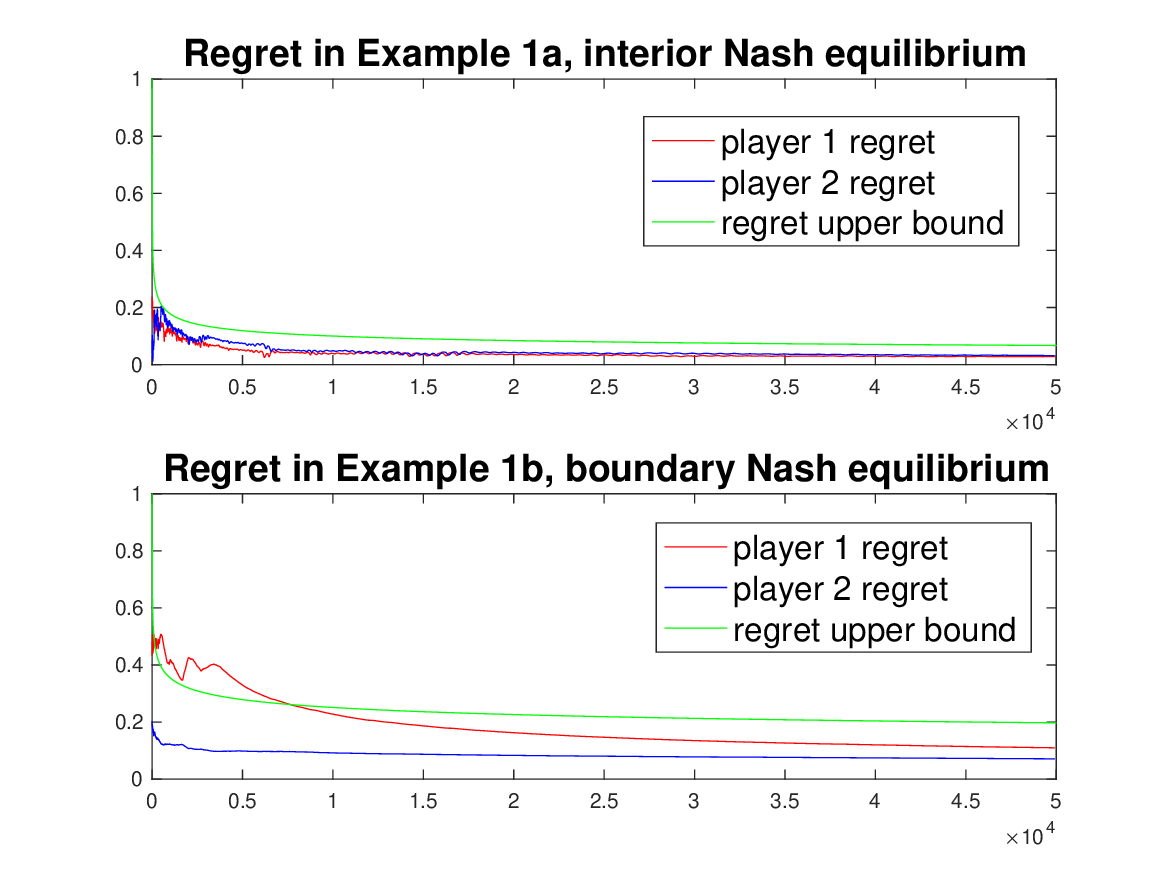}
\caption{Regret of the algorithm in Examples 1a) (top) and Example 1b.}
\label{fig:regret}
\end{subfigure}
\end{figure}

 {The regret of our algorithm for Examples 1a and 1b are shown in Figure \ref{fig:regret}. As expected from Corollary \ref{cor:regret_rate}, we observe higher regret  for  Example 1b, since in this case the least-norm Nash equilibrium is on the boundary. The regret rate derived from Corollary~\ref{cor:regret_rate} is shown, where $\varepsilon$ was set to  $0.04$.}

\subsection{Example 2} 
\begin{figure}[!t]
\centering
\begin{subfigure}
\centering
\includegraphics[width=\linewidth]{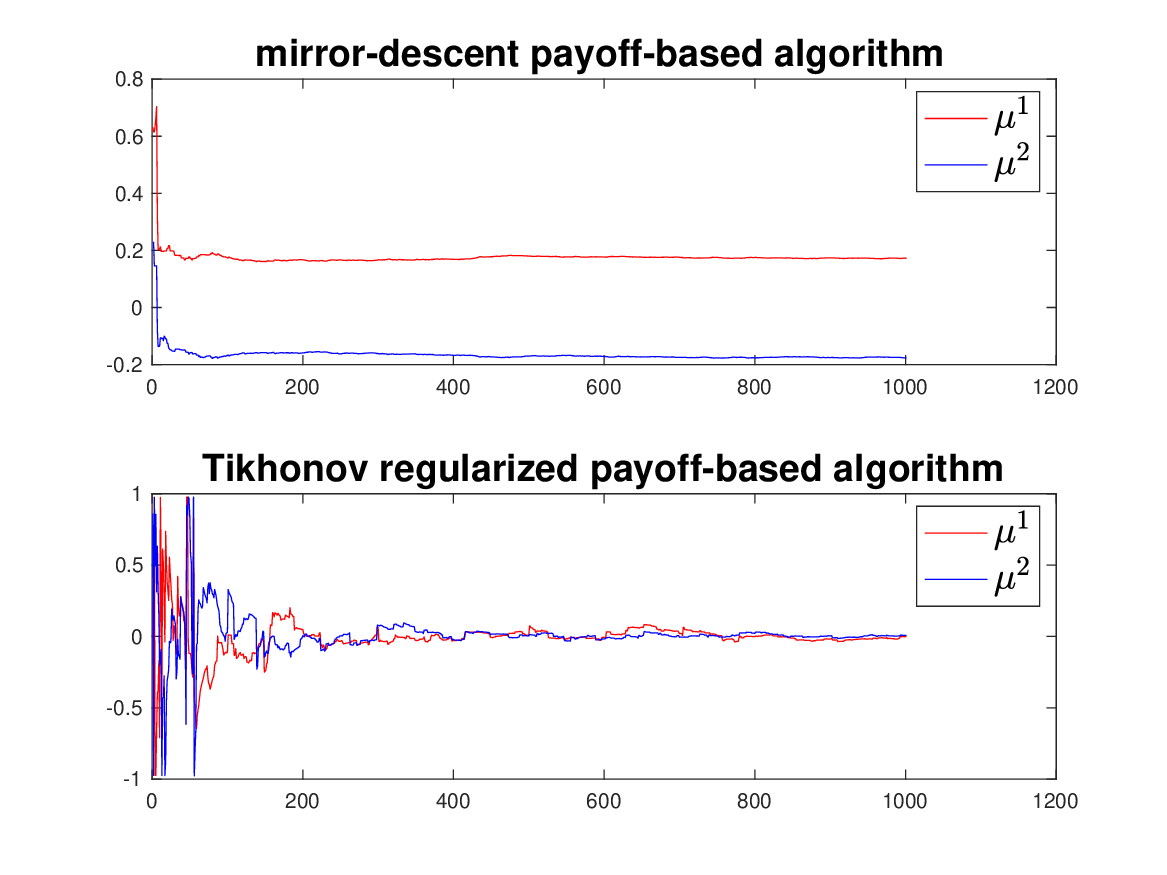}
\caption{The iterates in Example 2 for an arbitrary initial condition $\bmu(0)$.}
\label{fig:ex2_MD}
\end{subfigure}
\end{figure}

The iterates of the algorithm for Example \ref{example:nonunique} are provided in Fig. \ref{fig:ex2_MD}. The algorithm iterates converge to the least-norm Nash equilibrium. The iterates of \cite{bravo2018bandit} converge, but not to the least-norm Nash equilibrium. While the payoff-based mirror descent algorithm does not have guaranteed convergence for an arbitrary merely monotone game,  we can foresee the potential convergence of mirror descent in this particular example as follows.  Given the pseudo-gradient mapping $\mathbf{M}(a^1, a^2) = (a^1+a^2, a^1+a^2)^T$ in Example 2, we can directly verify that $\langle\mathbf{M}(\mathbf{a}) - \mathbf{M}(\mathbf{a}^*), \mathbf{a}-\mathbf{a}^*\rangle = (a^1-a^{1,*}+a^2-a^{2,*})^2$, where $\mathbf{a}^* = (a^{1,*}+a^{2,*}$ is a Nash equilibrium and, thus, $a^{1,*}+a^{2,*} = 0$. So, let us consider the distance between the iterates $\mathbf{a}(t+1) = \mbox{Proj}_{A}[\mathbf{a}(t) - \alpha_t\mathbf{M}(\mathbf{a}(t))]$ and $\mathbf{a}^* = \mbox{Proj}_{A}[\mathbf{a}^* - \alpha_t\mathbf{M}(\mathbf{a}^*)]$ in the first-order gradient procedure (a subclass of the mirror-descent methods): 
\begin{align*}
    &\|\mathbf{a}(t+1) - \mathbf{a}^*\|^2 \le \|\mathbf{a}(t) - \mathbf{a}^*\|^2 -2\alpha_t\langle\mathbf{M}(\mathbf{a}) - \mathbf{M}(\mathbf{a}^*), \mathbf{a}-\mathbf{a}^*\rangle + \alpha_t^2\|\mathbf{M}(\mathbf{a}) - \mathbf{M}(\mathbf{a}^*)\|\cr
    &\le (1+\alpha_t^2L^2)\|\mathbf{a}(t) - \mathbf{a}^*\|^2 - 2\alpha_t(a^1(t)-a^{1,*}+a^2(t)-a^{2,*})^2\cr
    &=(1+\alpha_t^2L^2-2\alpha_t)\|\mathbf{a}(t) - \mathbf{a}^*\|^2 - (4-\theta_t+\theta_t)\alpha_t(a^1(t)-a^{1,*})(a^2(t)-a^{2,*})\cr
    &\le (1+\alpha_t^2L^2-\alpha_t\theta_t/2 )\|\mathbf{a}(t) - \mathbf{a}^*\|^2 + O(\theta_t\alpha_t)
    \end{align*}
    for some $\theta_t\in(0,1)$ (analysis relevant constant) and where $L$ is the Lipschitz constant of $\mathbf{M}$. In the first inequality we used Lipshitz continuity of $\mathbf{M}$ and its property $\langle\mathbf{M}(\mathbf{a}) - \mathbf{M}(\mathbf{a}^*), \mathbf{a}-\mathbf{a}^*\rangle = (a^1-a^{1,*}+a^2-a^{2,*})^2$, whereas in the last inequality we used the relation $ -xy\le 1/2x^2+1/2y^2$ with $x = (a^1(t)-a^{1,*})$ and $y = (a^2(t)-a^{2,*})$. Thus, $\alpha_t$ and $\theta_t$ can be balanced in such a way that the Gladyshev's result (see Lemma 9 in Chapter 2 in \cite{polyak}) holds and one can conclude convergence of $\|\mathbf{a}(t+1) - \mathbf{a}^*\|^2$. We emphasize once again that as there are no guaranteed convergence of mirror descent for merely monotone cases, one needs to address this convergence case by case.

\begin{figure}[!t]
\centering
\begin{subfigure}
\centering
\includegraphics[width=\linewidth]{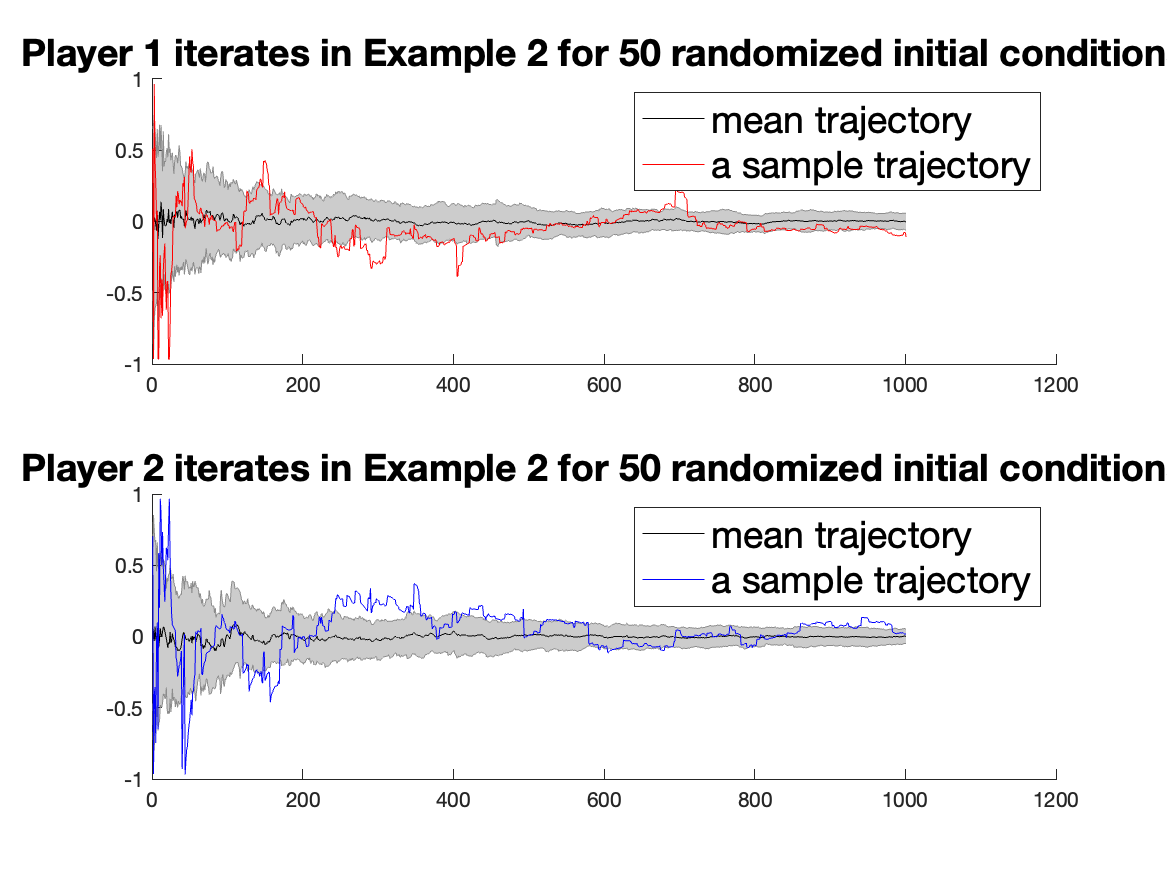}
\caption{Iterates in Example 2 - one standard deviation (in gray) of the 50 randomly selected initial condition is shown.}
\label{fig:randomized}
\end{subfigure}
\end{figure}

Figure~\ref{fig:randomized} illustrates the iterates of the proposed algorithm with 50 randomized initial conditions, $\bmu(0)$ for Example~\ref{example:nonunique}. The iterates converge to the least-norm Nash equilibrium as predicted by Theorem~\ref{th:main}.  {The mean of the iterates, the one-standard deviation, and one sample trajectory are illustrated in this figure.}

\subsection{Monotone game with varying dimension}\label{ex:varyingd}
\begin{figure}[!t]
\begin{subfigure}
		\centering
		\includegraphics[width=\linewidth]{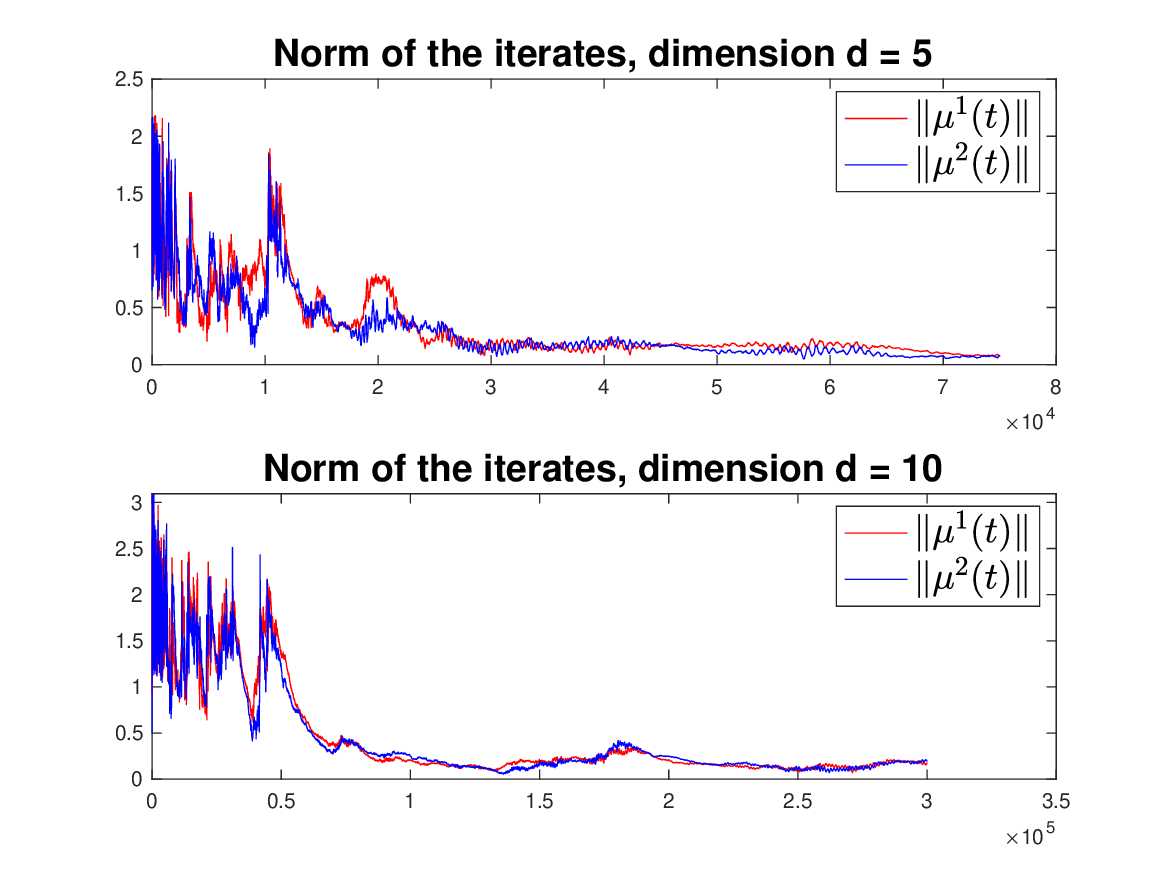}
		\caption{Algorithm \ref{alg:algorithm1} applied to the bilinear zero-sum matrix game with $J^1(\mathbf{a^1}, \mathbf{a^2}) = \mathbf{a^1}^T \mathbf{1}_{d} \mathbf{a^2}$, and $\Ab^i = [-1,1]^d$, $i=1,2$. }
	\label{fig:bilinear}
		\end{subfigure}
\end{figure}

{To illustrate the applicability of the proposed procedure to higher dimensional problems, we consider  a zero-sum game where the objective of player $1$ is $J^1(\mathbf{a^1}, \mathbf{a^2}) = \mathbf{a^1}^T \mathbf{1}_{d} \mathbf{a^2}$,  with $\mathbf{1}_{d} \in \R^{d\times d}$ being the matrix with unit 1 on all entries. The action sets are $\Ab^i =[-1,1]^d$. It can be verified that the game mapping is monotone and any pair of strategies satisfying $a^i \in \text{Null}(\mathbf{1}_d) \cap [-1,1]^d$, and $a^{-i}=\mathbf{0}_d$ is an equilibrium. Using Algorithm \ref{alg:algorithm1}, the iterates  converge to the least-norm Nash equilibrium, which is at $(\boldsymbol{0}, \boldsymbol{0}) \in \R^{2d}$. In Figure \ref{fig:bilinear}, we provide sample trajectories of the iterates of Algorithm \ref{alg:algorithm1}  for dimensions $d=5, 10$. }

\section{Conclusions}\label{sec:conclusions}
We considered the class of convex games with merely monotone pseudo-gradient and derived a payoff-based algorithm whose iterates converge to a least-norm Nash equilibrium. We furthermore established that our proposed algorithm is no-regret and quantified its regret rate by properly setting the algorithm parameters. An open future direction is developing alternative payoff-based algorithms with improved regret rates for this class of games and characterizing the regret rate as a function of the dimension of the action spaces.

\bibliographystyle{plain}
\bibliography{neurips_ref}

 \newpage

\begingroup
\let\clearpage\relax 
\onecolumn 
\endgroup


\appendix
%

 \section{Proof of Proposition~\ref{prop:sample_grad}}\label{app:sample_grad}
 First, let us derive the terms $\Qb^i$, $\Rb^i$ defined in Equation \eqref{eq:map_est}. Consider the gradient estimation in \eqref{eq:est_Gd}. Let  $p( \bx; \bmu, \sigma)=\prod_{i=1}^Np^i(x^i_1,\ldots,x^i_{d};\bmu^i,\sigma)$ denote the  density function of the joint distribution of the agents' sampling variables $(\bxi^1,\ldots,\bxi^N)$, where $p^i(x^i_1,\ldots,x^i_{d};\bmu^i,\sigma)$  is defined in \eqref{eq:density}. Given $\sigma > 0$, for any $i\in[N]$ define $ \tilde{J}^i : \R^{Nd} \rightarrow \R$ as
 \begin{align}
 	\label{eq:mixedJ}
 	\tilde{J}^i &(\bmu^1,\ldots,\bmu^N)= \int_{\mathbb R^{Nd}}J^i(\bx)p( \bx; \bmu, \sigma)d\bx.
 \end{align}
 Thus, $\tilde{J}^i$, $i\in[N]$, is the $i$-th player's cost function in mixed strategies, namely
 \[\tilde{J}^i (\bmu^1,\ldots,\bmu^N) = \E_{\bxi\sim \EuScript N(\bmu, \sigma)}J^i(\bxi).\]
 Let $\bmu(t) =(\bmu^1(t),\ldots,\bmu^N(t))$ where $\bmu^i(t) = (\mu^i_1,\ldots,\mu^i_d)$ is the current state of the agent $i$.  For $i\in[N]$, let us define $\tilde{\Mb}^i (\cdot)=(\tilde M^i_1(\cdot), \ldots, \tilde M^i_d(\cdot))^{\top}$
 as the $d$-dimensional mapping:
 \begin{align}\label{eq:mapp2}
 	\tilde{M}^i_{k} (\bmu,\sigma)=\frac{\partial {\tilde J^i(\bmu, \sigma)}}{\partial \mu^i_k}, \mbox{ for $k\in[d]$}.
 \end{align}
 This mapping is the  pseudo-gradient corresponding to the game with mixed strategy cost functions $\tilde{J}^i$ defined in \eqref{eq:mixedJ}.

 Next, observe that for all $i \in [N]$ we can equivalently write $\mb^i(t)$ defined in \eqref{eq:est_Gd} as
 \begin{align}
 	\label{eq:pbavmu}
 	\mb^i(t) = \Mb^i(\bmu(t))+\Qb^i(t) +\Rb^i(t),
 \end{align}
 where $\Qb^i$ and $\Rb^i$ are defined as
 \begin{align}
 	\label{eq:Qi}
 	&\Qb^i(t) =\underbrace{\tilde{\Mb}^i (\bmu(t)) -\Mb^i(\bmu(t))}_{\text{term 1}} + \underbrace{(J^i(\ab(t))-J^i(\bxi(t)))\frac{\bxi^i(t) -\bmu^i(t)}{\sigma^2_t}}_{\text{term 2}},\\
 	\label{eq:Ri}
 	&\Rb^i(t) = J^i(\bxi(t))\frac{\bxi^i(t) -\bmu^i(t)}{\sigma^2_t} - \tilde{\Mb}^i (\bmu(t)).
 \end{align}
 With the definitions of $\Qb^i$ and $\Rb^i$ in place, we proceed with establishing the properties of the stochastic terms $\Qb^i(t)$ and $\Rb^i(t)$ as stated in parts \ref{itm:R_mean}-\ref{itm:Q_norm} in the following three sections.

 \subsection{Part \ref{itm:R_mean} of Proposition~\ref{prop:sample_grad}}
 \begin{lem}
 	\label{lem:sample_grad}
 	Under Assumption\r\ref{assum:convex}
 	\begin{enumerate}
 		\item \label{itm:diff_int} $\tilde M^i_{k} (\bmu,\sigma) = \int_{\mathbb R^{Nd}}\frac{\partial {J^i(\bx)p( \bx; \bmu, \sigma)d\bx}}{\partial \mu^i_k}, \mbox{ for $k\in[d]$}$;
 		\item \label{itm:zero_mean} $\E\{J^i(\bxi^1(t),\ldots,\bxi^N(t))\frac{\xi^i_k(t) -\mu^i_k(t)}{\sigma^2_t}|\xi^i_k(t)\sim\EuScript N(\mu_k^i(t),\sigma_t)\}=\tilde M^i_{k} (\bmu(t),\sigma_t)$.
 	\end{enumerate}
 \end{lem}
 \begin{proof}
 	For Part ~\ref{itm:diff_int}, we  verify that  the differentiation under the integral sign in
 	\begin{align*}
 		\tilde M^i_{k} (\bmu,\sigma)=\frac{\partial {\int_{\mathbb R^{Nd}}J^i(\bx)p( \bx; \bmu, \sigma)d\bx}}{\partial \mu^i_k}, \mbox{ for $k\in[d]$},
 	\end{align*}
 	is justified. It can then readily be verified that  the statement~\ref{itm:zero_mean} holds, by taking the differentiation inside the integral. A sufficient condition for differentiation under the integral is that the integral of the formally differentiated function with respect to $\mu^i_k$  converges uniformly over the domain set of the parameter $\mu^i_k$, whereas the differentiated function is continuous (see \cite[Chapter 17]{zorich}). {Continuity of the functions $J^i(\bx)p( \bx; \bmu, \sigma)$ and
 		$\frac{\partial {J^i(\bx)p( \bx; \bmu, \sigma)}}{\partial \mu^i_k} = J^i(\bx)\frac{x^i_k -\mu^i_k}{\sigma^2}p( \bx; \bmu, \sigma)$
 		follows from continuity of $J^i(\bx)$ and $p( \bx; \bmu, \sigma)$. We demonstrate uniform convergence of the integral
 		$\int_{\R^{Nd}} J^i(\bx)\frac{x^i_k -\mu^i_k}{\sigma^2}p( \bx; \bmu, \sigma)d\bx$
 		over the set $\Ab$. As
 		\begin{align*}
 			&\int_{\R^{Nd}} J^i(\bx)\frac{x^i_k -\mu^i_k}{\sigma^2}p( \bx; \bmu, \sigma)d\bx  \\
 			&=\E\{J^i(\bxi)\frac{\xi^i_k -\mu^i_k}{\sigma^2}|\bxi\sim\EuScript N(\bmu,\sigma)\},
 		\end{align*}
 		we can apply the H\"older's inequality to conclude that
 		\begin{align}\label{eq:unifconv}
 			&\int_{\R^{Nd}} J^i(\bx)\frac{x^i_k -\mu^i_k}{\sigma^2}p( \bx; \bmu, \sigma)d\bx  \cr
 			&\le(\E\{(J^i)^2(\bxi)|\bxi\sim\EuScript N(\bmu,\sigma)\})^{\frac{1}{2}}(\E\{\frac{(\xi^i_k -\mu^i_k)^2}{\sigma^4}|\bxi\sim\EuScript N(\bmu,\sigma)\})^{\frac{1}{2}}\cr
 			& = \frac{1}{\sigma}(\E\{(J^i)^2(\bxi)|\bxi\sim\EuScript N(\bmu,\sigma)\})^{\frac{1}{2}}.
 		\end{align}
 		Next,
 		\begin{align*}
 			&\E\{{J^i}^2(\bxi)|\bxi\sim\EuScript N(\bmu,\sigma)\} = \int_{\R^{Nd}} {J^i}^2(\bx)p( \bx; \bmu, \sigma)d\bx\\
 			&=\int_{\mathbb R^{Nd}}{J^i}^2(\bx)\prod_{i=1}^{N} \frac{1}{(\sqrt{2\pi}\sigma)^{d}}\exp\left\{-\sum_{k=1}^{d}\frac{(x^i_k-\mu^i_k)^2}{2\sigma^2}\right\}d\bx\cr
 			& = 2^{Nd}\int_{\mathbb R^{Nd}}{J^i}^2(\bx)\prod_{i=1}^{N} \exp\left\{-\sum_{k=1}^{d}\frac{(x^i_k-\mu^i_k)^2}{4\sigma^2}\right\}\cr
 			&\qquad\qquad \times \frac{1}{(\sqrt{2\pi}2\sigma)^{Nd}}\prod_{i=1}^{N} \exp\left\{-\sum_{k=1}^{d}\frac{(x^i_k-\mu^i_k)^2}{4\sigma^2}\right\}d\bx.
 		\end{align*}
 		According to Assumption~\ref{assum:Lipschitz} Part ~\ref{itm:growth}, there exists the constants $c$  such that
 		\[{J^i}^2(\bx)\prod_{i=1}^{N} \exp\left\{-\sum_{k=1}^{d}\frac{(x^i_k-\mu^i_k)^2}{4\sigma^2}\right\}\le c\]
 		for any $\bx\in\R^{Nd}$ and $\bmu\in\Ab$.
 		Thus, taking into account that
 		\[\int_{\mathbb R^{Nd}}\frac{1}{(\sqrt{2\pi}2\sigma_{t})^{Nd}}\prod_{i=1}^{N} \exp\left\{-\sum_{k=1}^{d}\frac{(x^i_k-\mu^i_k)^2}{4\sigma^2_{t}}\right\}d\bx = 1,\]
 		we conclude the uniform convergence over $\Ab$ follows from~\eqref{eq:unifconv} and the fact that \[\int_{\mathbb R^{Nd}}{J^i}^2(\bx)\prod_{i=1}^{N} \frac{1}{(\sqrt{2\pi}\sigma)^{d}}\exp\left\{-\sum_{k=1}^{d}\frac{(x^i_k-\mu^i_k)^2}{2\sigma^2}\right\}d\bx\] converges uniformly over $\Ab^i$ according to Dirichlet's test for uniform convergence of integrals.}
 \end{proof}

 Note that, as a consequence of the above  lemma,  $\E\{\Rb^i(t)|\EuScript F_t\} = 0$ for all $i$ and any $t$ and thus, part~\ref{itm:R_mean} in Proposition~\ref{prop:sample_grad} is proven.  Furthermore, we can see that the term $\Rb^i(t)$ can be interpreted as a zero mean noise term in the gradient estimation. This implies that the expression $J^i(\bxi^1(t),\ldots,\bxi^N(t))\frac{\bxi^i(t) -\bmu^i(t)}{\sigma^2_t}$ is an unbiased one-point estimation of the gradient of the agent $i$'s cost function in mixed strategies, $\tilde{J}^i$, with respect to her decision variable.

 \subsection{Part \ref{itm:R_square} of Proposition~\ref{prop:sample_grad}}
 From the definition of $\Rb^i$ in \eqref{eq:Ri}, notice that
 \begin{align}
 	\label{eq:Rineq1}
 	\E\{\|\Rb^i(t)\|^2|\EuScript F_t\}&=\E\{\|\Rb^i(t)\|^2|\bxi(t)\sim \EuScript(\bmu(t),\sigma_t)\}\cr
 	&\le\sum_{k=1}^{d}\int_{\mathbb R^{Nd}}{J^i}^2(\bx)\frac{(x^i_k - \mu^i_k(t))^2}{\sigma^4_t} p(\bx;\bmu(t),\s_t)d\bx
 \end{align}
 Next, according to the H\"older's inequality ($\E|XY|\le (\E(X^2))^{1/2}(\E(Y^2))^{1/2})$, we have
 \begin{align}\label{eq:boundR0}
 	\int_{\mathbb R^{Nd}}&{J^i}^2(\bx)\frac{(x^i_k - \mu^i_k(t))^2}{\sigma^4_t} p(\bx;\bmu(t),\s_t)d\bx \cr
 	&= \E\{{J^i}^2(\bxi(t))\frac{(\xi^i_k - \mu^i_k(t))^2}{\sigma^4_t}\}\le (\E\{{J^i}^4(\bxi(t))\})^{1/2}\left(\E\{\frac{(\xi^i_k - \mu^i_k(t))^4}{\sigma^8_t}\}\right)^{1/2}\cr
 	& = \frac{1}{\sigma_t^2}(\E\{{J^i}^4(\bxi(t))\})^{1/2},
 \end{align}
 where in the inequality above for the safe of notation simplicity we used $\E\{\cdot\} = \E\{\cdot |\bxi(t)\sim \EuScript(\bmu(t),\sigma_t)\}$.
 We continue analysing the term $\E\{{J^i}^4(\bxi(t))\}$:
 \begin{align}\label{eq:boundR}
 	\E\{{J^i}^4(\bxi(t))\} = &\int_{\mathbb R^{Nd}}{J^i}^4(\bx) p(\bx;\bmu(t),\s_t)d\bx\cr
 	& =\int_{\mathbb R^{Nd}}{J^i}^4(\bx)\prod_{k=1}^{N} \frac{1}{(\sqrt{2\pi}\sigma_{t})^{d}}\exp\left\{-\sum_{k=1}^{d}\frac{(x^i_k-\mu^i_k(t))^2}{2\sigma^2_{t}}\right\}d\bx\cr
 	& = 2^{Nd}\int_{\mathbb R^{Nd}}{J^i}^4(\bx)\prod_{k=1}^{N} \exp\left\{-\sum_{k=1}^{d}\frac{(x^i_k-\mu^i_k(t))^2}{4\sigma^2_{t}}\right\}\cr
 	&\qquad\qquad \times \frac{1}{(\sqrt{2\pi}2\sigma_{t})^{Nd}}\prod_{k=1}^{N} \exp\left\{-\sum_{k=1}^{d}\frac{(x^i_k-\mu^i_k(t))^2}{4\sigma^2_{t}}\right\}d\bx.
 \end{align}
 According to Assumption~\ref{assum:Lipschitz} Part ~\ref{itm:growth}, there exists the constants $C$  such that
 \[{J^i}^4(\bx)\prod_{k=1}^{N} \exp\left\{-\sum_{k=1}^{d}\frac{(x^i_k-\mu^i_k)^2}{4\sigma^2_{t}}\right\}\le C\]
 for any $t$, $\bx\in\R^{Nd}$, and $\bmu\in\Ab$.
 Thus, taking into account that
 \[\int_{\mathbb R^{Nd}}\frac{1}{(\sqrt{2\pi}2\sigma_{t})^{Nd}}\prod_{k=1}^{N} \exp\left\{-\sum_{k=1}^{d}\frac{(x^i_k-\mu^i_k(t))^2}{4\sigma^2_{t}}\right\}d\bx = 1,\]
 we conclude that  $\E\{{J^i}^4(\bxi(t))\}\le 2^{Nd}C$ and the result follows from \eqref{eq:boundR0}.

 \subsection{Part \ref{itm:Q_norm} of Proposition~\ref{prop:sample_grad}}
 Recall the definition of $\Qb^i$ in \eqref{eq:est_Gd}.  We address proving this part of the proposition by first bounding the last term $(J^i(\ab(t))-J^i(\bxi(t)))\frac{\bxi^i(t) -\bmu^i(t)}{\sigma^2_t}$ in the biased bias term $\Qb^i(t)$.  Our next lemma bounds the norm of this bias term in terms of $\rho_t$ and $\sigma_t$.

 \begin{lem}\label{lem:projTerm}
 	Under Assumptions ~\ref{assum:convex} and ~\ref{assum:Lipschitz} we have
 	\begin{align}
 		\E\left\{|J^i(\ab(t))-J^i(\bxi(t))|\frac{\|\bxi^i(t) -\bmu^i(t)\|}{\sigma^2_t}\;|\EuScript F_t\right\}=O\left(\left(\frac{e^{-\frac{{\rho}_t^2}{2\sigma_t^2}}}{\sigma_t^{Nd+1}}\right)^{\frac{1}{2}}\right) \mbox{ almost surely. }
 	\end{align}
 \end{lem}
 \begin{proof}
 	As $\ab(t) = \Proj_{\Ab}\bxi(t)$, we conclude that almost surely
 	\begin{align}\label{eq:projterm0}
 		&\E\left\{|J^i(\ab(t))-J^i(\bxi(t))|\frac{\|\bxi^i(t) -\bmu^i(t)\|}{\sigma^2_t}\;|\EuScript F_t\right\} \cr
 		&= \E\left\{|J^i(\ab(t))-J^i(\bxi(t))|\frac{\|\bxi^i(t) -\bmu^i(t)\|}{\sigma^2_t}| \bxi(t)\sim \EuScript N(\bmu(t),\s_t) \right\}\cr
 		&\le (\E\left\{|J^i(\ab(t))-J^i(\bxi(t))|^2\right\})^{\frac{1}{2}}\left(\E\left\{\frac{\|\bxi(t)-\bmu(t)\|^2}{\sigma^4_t}\right\}\right)^{\frac{1}{2}}\cr
 		& =\frac{\sqrt{Nd}}{\sigma_t}(\E\left\{|J^i(\ab(t))-J^i(\bxi(t))|^2\right\})^{\frac{1}{2}},
 	\end{align}
 	where, as before, we used the H\"older's inequality ($\E|XY|\le (\E(X^2))^{1/2}(\E(Y^2))^{1/2})$ and the notation $\E\{\cdot\} = \E\{\cdot |\bxi(t)\sim \EuScript(\bmu(t),\sigma_t)\}$.
 	Next, we estimate the term $\E\left\{|J^i(\ab(t))-J^i(\bxi(t))|^2\right\}$ as follows.
 	\begin{align}\label{eq:projterm}
 		&\E\left\{|J^i(\ab(t))-J^i(\bxi(t))|^2\right\} \cr
 		&=\int_{\Ab}|J^i(\Proj_{\Ab}\bx)-J^i(\bx)|^2p(\bx;\bmu(t),\sigma_t)d\bx + \int_{\R^{Nd}\setminus \Ab}|J^i(\Proj_{\Ab}\bx)-J^i(\bx)|^2p(\bx;\bmu(t),\sigma_t)d\bx\cr
 		&=\int_{\R^{Nd}\setminus \Ab}|J^i(\Proj_{\Ab}\bx)-J^i(\bx)|^2p(\bx;\bmu(t),\sigma_t)d\bx\cr
 		&\le 2\int_{\R^{Nd}\setminus \Ab}(J^i(\Proj_{\Ab}\bx))^2p(\bx;\bmu(t),\sigma_t)d\bx + 2\int_{\R^{Nd}\setminus \Ab}(J^i(\bx))^2p(\bx;\bmu(t),\sigma_t)d\bx\cr
 		&\le 2K\Pr\{\bxi(t)\in\R^{Nd}\setminus \Ab\}) + 2\int_{\R^{Nd}\setminus \Ab}(J^i(\bx))^2p(\bx;\bmu(t),\sigma_t)d\bx,
 	\end{align}
 	where the second equality above is due to the fact that $\Proj_{\Ab}\bx=\bx$ for any $\bx\in \Ab$. The first inequality was obtained by taking into account that $|J^i(\Proj_{\Ab}\bx)-J^i(\bx)|^2\le 2 (J^i(\Proj_{\Ab}\bx))^2 + 2(J^i(\bx))^2$, whereas the last one is due to the inequality $(J^i(\Proj_{\Ab}\bx))^2\le K$ for some constant $K$ and any $\bx$ (since $\Ab$ is compact and $J^i$ is continuous).

 	Thus, let us estimate $\Pr\{\bxi(t)\in\R^{Nd}\setminus \Ab\}$ above. The idea is that since $\bxi(t)$ is sampled from the Gaussian distribution with mean $\bmu(t)$, $\bxi(t)$ concentrates around its mean $\bmu(t)$ with high probability. Since the mean is projected onto a shrunk version of the set $\Ab$, namely, $(1-\rho_t) \Ab$, by appropriately tuning $\rho_t$ and the variance of the distribution $\sigma_t$ we can ensure that $\bxi(t)$ stays within the original feasible set with high probability.
 	
 	Let $\EuScript O_{\rho_t}(\xb) = \{\boldsymbol y\in\R^{Nd}| \|\boldsymbol y-\xb\|^2<{\rho}_t^2\}$ denote the $\rho_t$-neighborhood of the point $\xb\in\Ab$. Hence, $\sup_{\boldsymbol y\notin \EuScript O_{\rho_t}(\xb)}-\|\boldsymbol y-\xb\|^2=-{\rho}_t^2$ . Then, taking into account the fact that $\EuScript O_{\rho_t}(\xb)$ is contained in $\Ab$ and $\rho_t<1$, we obtain that for any $t$ and any bounded $\sigma>\sigma_t$:
 	\begin{align}\label{eq:probterm}
 		\Pr&\{\bxi(t)\in\R^{Nd}\setminus \Ab\}\le\Pr\{\bxi(t)\in\R^{Nd}\setminus \EuScript O_{\rho_t}(\bmu(t))\}\cr
 		& = \int_{\boldsymbol y\notin \EuScript O_{\rho_t}(\bmu(t))}\frac{1}{(2\pi)^{Nd/2}\sigma_t^{Nd}}
 		\exp\left\{-\frac{\|\boldsymbol y-\bmu(t)\|^2}{2\sigma_t^2}\right\}d\boldsymbol y \cr
 		&=\int_{\boldsymbol y\notin \EuScript O_{\rho_t}(\bmu(t))}\exp\left\{-\|\boldsymbol y-\bmu(t)\|^2\left(\frac{1}{2\sigma_t^2} - \frac{1}{2\sigma^2}\right)\right\}\cr
 		&\qquad\qquad\times\frac{\sigma^{Nd}}{\sigma_t^{Nd}}\frac{1}{(2\pi)^{{Nd}/2}\sigma^{Nd}}\exp\left\{-\frac{\|\boldsymbol y-\bmu(t)\|^2}{2\sigma^2}\right\}d\boldsymbol y\cr
 		&\le \exp\left\{-\rho^2_t\left(\frac{1}{2\sigma_t^2} - \frac{1}{2\sigma^2}\right)\right\}\frac{\sigma^{Nd}}{\sigma_t^{Nd}}\cr
 		&\quad\times\int_{\boldsymbol y\notin \EuScript O_{\rho_t}(\bmu(t))}\frac{1}{(2\pi)^{Nd/2}\sigma^{Nd}}\exp\left\{-  \frac{\|\boldsymbol y-\bmu(t)\|^2}{2\sigma^2}\right\}d\boldsymbol y\cr
 		&\le k_1 \frac{e^{-\frac{{\rho}_t^2}{2\sigma_t^2}}}{\sigma_t^{Nd}},
 	\end{align}
 	for some finite $k_1>0$. The last inequality holds because
 	\[\int_{\boldsymbol y\notin \EuScript O_{\rho_t}(\bmu(t))}\frac{1}{(2\pi)^{{Nd}/2}\sigma^{Nd}}\exp\left\{-  \frac{\|\boldsymbol y-\bmu(t)\|^2}{2\sigma^2}\right\}d\boldsymbol y\le 1\]
 	and, thus, due to the diminishing $\rho_t$ there exists $0<k_1<\infty$:
 	\[\int_{\boldsymbol y\notin \EuScript O_{\rho_t}(\bmu(t))}\frac{e^{ \frac{{\rho}_t^2}{2\sigma^2}}\sigma^{Nd}}{(2\pi)^{N/2}\sigma^{Nd}}\exp\left\{-  \frac{\|\boldsymbol y-\bmu(t)\|^2}{2\sigma^2}\right\}d\boldsymbol y\le k_1.\]
 	
 	Taking into account Assumption~\ref{assum:Lipschitz}~\ref{itm:growth}, we conclude existence of some constant $C_1$ such that $(J^i(\bx))^2\exp\left\{-\frac{\|\boldsymbol x-\bmu(t)\|^2 }{4\sigma_t^2}\right\}\le C_1$ for any $\bx$. Thus,
 	\begin{align*}
 		&\int_{\R^{Nd}\setminus \Ab}(J^i(\bx))^2p(\bx;\bmu(t),\sigma_t)d\bx\cr
 		&\le C_1\int_{\boldsymbol y\in \R^{Nd}\setminus \Ab}\frac{1}{(2\pi)^{Nd/2}\sigma_t^{Nd}}
 		\exp\left\{-\frac{\|\boldsymbol x-\bmu(t)\|^2 }{4\sigma_t^2}\right\}d\boldsymbol x\cr
 		& \le C_1\int_{\boldsymbol x\notin \EuScript O_{\rho_t}(\bmu(t))}\frac{1}{(2\pi)^{Nd/2}\sigma_t^{Nd}}
 		\exp\left\{-\frac{\|\boldsymbol x-\bmu(t)\|^2 }{4\sigma_t^2}\right\}d\boldsymbol x.
 	\end{align*}
 	Hence, analogously to \eqref{eq:probterm}, we have
 	\begin{align}\label{eq:probterm1}
 		&\int_{\R^{Nd}\setminus \Ab}(J^i(\bx))^2p(\bx;\bmu(t),\sigma_t)d\bx \cr
 		& \le C_1\int_{\boldsymbol x\notin \EuScript O_{\rho_t}(\bmu(t))}\frac{1}{(2\pi)^{Nd/2}\sigma_t^{Nd}}
 		\exp\left\{-\frac{\|\boldsymbol x-\bmu(t)\|^2 }{4\sigma_t^2}\right\}d\boldsymbol x\cr
 		&=\int_{\boldsymbol x\notin \EuScript O_{\rho_t}(\bmu(t))}\exp\left\{-\|\boldsymbol x-\bmu(t)\|^2\left(\frac{1}{4\sigma_t^2} - \frac{1}{4\sigma^2}\right)\right\}\cr
 		&\qquad\qquad\times\frac{(\sqrt{2}\sigma)^{Nd}}{\sigma_t^{Nd}}\frac{1}{(2\pi)^{{Nd}/2}(\sqrt{2}\sigma)^{Nd}}\exp\left\{-\frac{\|\boldsymbol x-\bmu(t)\|^2}{4\sigma^2}\right\}d\boldsymbol x\cr
 		&\le \exp\left\{-\rho^2_t\left(\frac{1}{4\sigma_t^2} - \frac{1}{4\sigma^2}\right)\right\}\frac{(\sqrt{2}\sigma)^{Nd}}{\sigma_t^{Nd}}\cr
 		&\quad\times\int_{\boldsymbol x\notin \EuScript O_{\rho_t}(\bmu(t))}\frac{1}{(2\pi)^{Nd/2}(\sqrt{2}\sigma)^{Nd}}\exp\left\{-  \frac{\|\boldsymbol x-\bmu(t)\|^2}{4\sigma^2}\right\}d\boldsymbol x\cr
 		&\le k_2 \frac{e^{-\frac{{\rho}_t^2}{2\sigma_t^2}}}{\sigma_t^{Nd}}= O\left(\frac{e^{-\frac{{\rho}_t^2}{2\sigma_t^2}}}{\sigma_t^{Nd}}\right),
 	\end{align}
 	where the last inequality is due to the fact that
 	\begin{align*}
 		&\int_{\boldsymbol x\notin \EuScript O_{\rho_t}(\bmu(t))}\frac{1}{(2\pi)^{Nd/2}(\sqrt{2}\sigma)^{Nd}}\exp\left\{-  \frac{\|\boldsymbol x-\bmu(t)\|^2}{4\sigma^2}\right\}d\boldsymbol x\\
 		&\le \int_{\R^{Nd}}\frac{1}{(2\pi)^{Nd/2}(\sqrt{2}\sigma)^{Nd}}\exp\left\{-  \frac{\|\boldsymbol x-\bmu(t)\|^2}{4\sigma^2}\right\}d\boldsymbol x=1.
 	\end{align*}
 \end{proof}

 From Lemma~\ref{lem:projTerm}, it follows that the second term of $\Qb^i$ defined in \eqref{eq:Qi} converges to zero exponentially fast as $\lim_{t\to\infty}\frac{\rho_t}{\s_t} = \infty$. We proceed with bounding the first term in \eqref{eq:Qi}.

 Taking into account Assumption~\ref{assum:Lipschitz} Part\ref{itm:lip} and the mean value theorem, we have
 \begin{align}\label{eq:Qterm1}
 	\|\tilde{\Mb}^i (\bmu(t)) -\Mb^i(\bmu(t))\|&=\|\int_{\R^{Nd}}[\Mb^i(\bx)-\Mb^i(\bmu(t))]p(\bx;\bmu(t),\sigma_t)d\bx\| \cr
 	&\le\int_{\R^{Nd}}\|\Mb^i(\bx) - \Mb^i(\bmu(t))\| p(\bx;\bmu(t),\sigma_t) d\bx\cr
 	&=\int_{\R^{Nd}}\|(\nabla\Mb^i(\tilde{\bx}),\bx - \bmu(t))\| p(\bx;\bmu(t),\sigma_t) d\bx\cr
 	&\le\int_{\R^{Nd}}\|\nabla\Mb^i(\tilde{\bx})\|\|\bx - \bmu(t)\| p(\bx;\bmu(t),\sigma_t) d\bx,
 \end{align}
 where $\nabla\Mb^i(\tilde{\bx})$ is the gradient of the vector function $\Mb^i$ at the point $\tilde{\bx} = \bx+\theta(\bmu(t)-\bx)$ with $\theta\in[0,1]$ being a constant dependent on $\bx$ and $\bmu(t)$.
 Next, applying the H\"older's inequality to~\eqref{eq:Qterm1} and using the notations $\E\{\cdot\} = \E\{\cdot |\bxi(t)\sim \EuScript(\bmu(t),\sigma_t)\}$ and $\tilde{\bxi} = \bxi+\theta(\bmu(t)-\bxi)$,  we obtain
 \begin{align}\label{eq:Qterm2}
 	\|\tilde{\Mb}^i (\bmu(t)) -\Mb^i(\bmu(t))\|&\le \E\{\|\nabla\Mb^i(\tilde{\bxi})\|\|\bxi - \bmu(t)\|\} \cr
 	&\le (\E\{\|\nabla\Mb^i(\tilde{\bxi})\|^2\})^{\frac{1}{2}}(\E\{\|\bxi - \bmu(t)\|^2\})^{\frac{1}{2}}\cr
 	& = \sqrt{Nd}\sigma_t (\E\{\|\nabla\Mb^i(\tilde{\bxi})\|^2\})^{\frac{1}{2}}.
 \end{align}
 where the last equality is due to the fact that the first central absolute moment of a random variable with a normal distribution $\EuScript N(\mu,\sigma)$ is $O(\sigma)$.
 Taking into account Assumption~\ref{assum:Lipschitz} Part \ref{itm:growth} we can conclude that $\|\nabla\Mb^i({\bx})\|^2=O(\exp\{\|\bx\|^{\alpha'}\})$ as $\|\bx\|\to\infty$, where $\alpha'<2$. Hence, analogously to~\eqref{eq:boundR}, we can demonstrate existence of some constant $C_1$ such that  $\E\{\|\nabla\Mb^i(\tilde{\bxi})\|^2\}\le C_1$. Thus, \eqref{eq:Qterm2} implies the result.

 \subsection{Uniform Distribution for Gradient Estimation} \label{appendix:uniform}
 Let the  agent $i$ sample the random vector $\bxi^i(t)$ according to the multidimensional uniform distribution over the unit sphere.
 Then the action the agent $i$ chooses is $\ab^i(t) = \bmu^i(t)+\sigma_t\bxi^i(t)$, where $\sigma_t$ is some time dependent parameter. As before, the cost value $ J^i(t)$ at the joint action $\ab(t)=(\ab^1(t),\ldots,\ab^N(t))\in \Ab$, $ J^i(t) =J^i(\ab(t)) $ is revealed to each agent $i$. Based on this information, the agent $i$ estimates the local gradient $\frac{\partial J^i}{\partial \bmu^i}$ at the point of the joint state $\bmu(t)=(\bmu^1(t),\ldots,\bmu^N(t))$ as:
 \begin{align}\label{eq:est_ud}
 	\mb^i(t) = \frac{d}{\s_t} J^i(t)\bxi^i(t).
 \end{align}

 We note that ensuring the feasibility of actions as well as the properties of the terms $\Qb^i$ and $\Rb^i$ in the case of the estimations based on the uniform distribution in \eqref{eq:est_ud} are investigated in the works of \cite{flaxman2005online}, \cite{bravo2018bandit}.  In particular, as shown in \cite[Section 4]{bravo2018bandit} the gradient estimator has the same order of variance and bias as those derived in Proposition~\ref{prop:sample_grad}.

 \section{Proof of Lemma~\ref{lem:step_constants}}\label{app:step}
 \begin{proof}
 	The series $\sum_{t = 0}^\infty 1/t^p$ converges for $p > 1$ and diverges otherwise. Thus, it is trivial to verify all statements for Assumption~\ref{assum:parameter2}, except Part \ref{itm:omega}. To show the condition for that part, let us consider the term $(\e_t-\e_{t-1})^2$ in the first summand of Part \ref{itm:omega},  namely,  $\sum_{t=0}^{\infty}\frac{(\e^{1+\varepsilon}_t-\e^{1+\varepsilon}_{t-1})^2}{\e^{3}_t\g_t}$:
 	\begin{align*}
 		(\e_t-\e_{t-1})^2 &= (t^{-e} - (t-1)^{-e})^2\quad \text{(multiply by $\frac{t^{2e}}{t^{2e}}$)} \\
 		&= {\big((1-1/t)^{-e} -1\big)^2}/{t^{2e}}\quad\quad \text{(do Taylor  approximation)}\\
 		& = {( 1+ e/t + O(t^{-2}) -1)^2}/{t^{2e}} = O(t^{-2-2e}).
 	\end{align*}
 	Combining the above with the denominator $\e^{3}_t\g_t$, we obtain that $\sum_{t=0}^{\infty}\frac{(\e_t-\e_{t-1})^2}{\e^{3}_t\g_t}$ converges if $g+e<1$. Repeating the same analysis for $\frac{(\rho_t-\rho_{t-1})^2}{\e^{5+2\varepsilon}_t\g_t}$, we obtain $\sum_t \frac{(\rho_t-\rho_{t-1})^2}{\e^{5+2\varepsilon}_t\g_t}$ converges if $2+2r-5e-2\varepsilon e-g > 1$.  Note, however, that since $\varepsilon$ can be chosen arbitrary small, the inequality $2+2r-5e-g > 1$ implies the initial $2+2r-5e-2\varepsilon e-g > 1$ for sufficiently small $\varepsilon>0$. Thus, Part \ref{itm:omega} is verified.
 \end{proof}

\section{Proof of Proposition~\ref{th:Tikhonov}}
 \label{app:y_proof}
We use the following well-established result, see, for example, \cite{FaccPang1}.
 \begin{lem}\label{lem:VI_sol}
 	Consider a mapping $\boldsymbol T(\cdot): \R^d \to \R^d$ and a convex closed set $Y \subseteq \R^d$. Given $\theta > 0$,
 	\begin{align}
 		\label{ex:projection_VI}
 		\by^* \in \text{SOL}(Y, \boldsymbol T) \iff \by^* = \Proj_Y(\by^* - \theta \boldsymbol T(\by^*)).
 	\end{align}
 \end{lem}

 \begin{proof}

 {The proof is motivated by the work on solving monotone stochastic variational inequalities \cite{Koshal}. In contrast to the above work, we have an additional regularization term $\rho_t$. This regularization arises because we need to estimate the gradient using zero-order information at feasible points of the action set. }

 Let $\ba^*$ be the least-norm solution of $VI(\Ab,\Mb)$. Moreover, let $\ba^p$ be the projection of $\ba^*$ onto the set $(1-\rho_t)\Ab$. Next, let $\by(t)$ be the unique solution of the doubly regularized inequality, namely $\by(t)\in SOL((1-\rho_t)\Ab,\Mb + \e_t \mathbf{I}_{Nd})$. Thus, we conclude that
 \[\la\Mb(\ba^*), \by(t) - \ba^*\ra\ge 0,\]
 \[\la\Mb(\by(t)) + \e_t\by(t), \ba^p - \by(t)\ra\ge 0.\]
 Using the above two inequalities we obtain that
 \begin{align*}
 	0  \le &\la\Mb(\ba^*), \by(t) - \ba^*\ra +  \la\Mb(\by(t)) + \e_t\by(t), \ba^* - \by(t)\ra +  \la\Mb(\by(t)) + \e_t\by(t), \ba^p - \ba^*\ra\\
 	=&-\la\Mb(\ba^*)-\Mb(\by(t)), \ba^* - \by(t)\ra + \e_t\la\by(t),\ba^* - \by(t)\ra+\la\Mb(\by(t)) + \e_t\by(t), \ba^p - \ba^*\ra\\
 	\le&\e_t\la\by(t),\ba^*\ra - \e_t\|\by(t)\|^2+\la\Mb(\by(t)), \ba^p - \ba^*\ra + \e_t\la\by(t), \ba^p - \ba^*\ra,
 \end{align*}
 where in the last inequality we took into account monotonicity of $\Mb$, Hence,
 \begin{align*}
 	\e_t\|\by(t)\|^2 & \le \e_t(\by(t),\ba^*)+\la\Mb(\by(t)), \ba^p - \ba^*\ra + \e_t\la\by(t), \ba^p - \ba^*\ra \\
 	& \le \e_t\|\by(t)\|\|\ba^*\|+l\| \ba^p - \ba^*\| + \e_t\|\by(t)\| \|\ba^p - \ba^*\| \\
 	& =  \e_t\|\by(t)\|\|\ba^*\|+l O(\rho_t) + \e_t\|\by(t)\| O(\rho_t),
 \end{align*}
 where in the first inequality we used Remark~\ref{rem:LipschJ} and in the second one we applied Lemma~\ref{lem:A_r}.
 Hence,
 \[\|\by(t)\|^2  \le \|\by(t)\|\|\ba^*\|+l O\left(\frac{\rho_t}{\e_t}\right) + \|\by(t)\| O(\rho_t).\]
 By taking the upper limit as $t\to\infty$  above and due to the settings for $\e_t$ and $\rho_t$, we obtain
 \begin{align*}\overline{\lim}_{t\to\infty}[\|\by(t)\|^2  ]\le &\|\ba^*\|\overline{\lim}_{t\to\infty}\|\by(t)\|+ \overline{\lim}_{t\to\infty}\|\by(t)\|\overline{\lim}_{t\to\infty} O(\rho_t).
 \end{align*}
 It implies that $\overline{\lim}_{t\to\infty}\|\by(t)\|  \le \|\ba^*\|$,
 and, thus, the sequence $\|\by(t)\|$ is upper bounded. Moreover, any accumulation point of this sequence is bounded above by the Euclidean
 norm of $\ab^*$, the least-norm solution of $VI(\Ab,\Mb)$. Further, according to the fact that
 \[\by(t)=\Proj_{(1-\rho_t)\Ab}[\by(t)-\theta(\Mb + \e_t \by(t))], \mbox{ for any } \theta>0,\]
 and that the function $\Proj_{(1-\rho)\Ab}(\xb)$ is continuous\footnote{It follows from the fact that for any $\rho$, $\rho'>0$ and $\xb_1$, $\xb_2\in\R^{Nd}$ we get
 	$
 	\|\Proj_{(1-\rho)\Ab}[\xb_1]-\Proj_{(1-\rho')\Ab}[\xb_2]\|
 	\le \|\Proj_{(1-\rho)\Ab}[\xb_1]-\Proj_{(1-\rho)\Ab}[\xb_2]\| + O(|\rho-\rho'|)
 	\le \|\xb_1-\xb_2\| + O(|\rho-\rho'|).
 	$} in both $\rho$ and $\xb$  any accumulation point of the sequence $\by(t)$ is a solution of $VI(\Ab,\Mb)$.
 Hence, as $\|\by(t)\|\le\|\ab^*\|$, where $\ab^*$ is the least-norm solution of $VI(\Ab,\Mb)$, we conclude that $\by(t)$ converges to such least-norm solution of $VI(\Ab,\Mb)$.
 \end{proof}
 
\section{Proof of Proposition~\ref{prop:zt}}\label{app:zt}
 
Our goal is to establish $\zb(t)$ defined in~\eqref{eq:crmd} converges to $\by(t)$. Then, from Proposition~\ref{th:Tikhonov} we can establish convergence to a Nash equilibrium for the sequence $\zb(t)$. 

 First, let us focus on proving the above lemma by establishing some useful properties of projecting onto the sets $(1-\rho_t)\Ab$.
 \begin{lem}\label{lem:A_r}
 For any $\xb\in\R^{Nd}$ the following holds:
 \[\|\Proj_{(1-\rho_{t-1})\Ab}\xb - \Proj_{(1-\rho_{t})\Ab}\xb\| = O(|\rho_{t-1}-\rho_t|).\]
 \end{lem}
 \begin{proof}
 Without loss of generality,  assume  $\xb\notin (1-\rho_{t-1})\Ab$ (otherwise, $\|\Proj_{(1-\rho_{t-1})\Ab}\xb - \Proj_{(1-\rho_{t})\Ab}\xb\| = 0$).
 Due to convexity of the set $\Ab\subset \R^{Nd}$ there exists a convex function $g:\R^{Nd}\to\R$ such that $\Ab = \{\xb: g(\xb)\le 0\}$, so that $(1-\rho_t)\Ab = \{\xb: g(\xb)\le -\rho_t\}$ for any $t$. Moreover, define $\xb': = \Proj_{(1-\rho_{t})\Ab}\xb,$ and observe that $ \Proj_{(1-\rho_{t-1})\Ab}\xb = \Proj_{(1-\rho_{t-1})\Ab}\xb'$. Thus, we have $\|\Proj_{(1-\rho_{t-1})\Ab}\xb - \Proj_{(1-\rho_{t})\Ab}\xb\|= d,$ where $d$ is the optimal cost corresponding to the following optimization problem
 \begin{align*}
 	d := \min_{\yb} &\|\yb - \xb'\|, \\
 	\mbox{s.t. } &g(\yb) = -\rho_{t-1}.
 \end{align*}
 The optimization problem has a solution $\yb^*$ for which the gradient of the corresponding Lagrangian is zero, namely
 \[\frac{(\yb^* - \xb')}{\|\yb^* - \xb'\|} + \lambda \nabla g(\yb^*) = \boldsymbol 0,\]
 where $\lambda>0$ is the dual multiplier of the problem under consideration. Due to Assumption~\ref{assum:convex} and the choice of $\rho_t$, that guarantees nonempty interior of $(1-\rho_t)\Ab$ for all $t$, the Slator's condition for the constraints $g(\xb)\le -\rho_t$ holds for all $t$. Hence, for any $\xb\in\R^{Nd}$ there exists a constant $\Lambda>0$ such that $\lambda<\Lambda$ (see \cite{ConvOpt}).
 Thus, we conclude that
 \[\nabla g(\yb^*) = -\frac{(\yb^* - \xb')}{ \lambda \|\yb^* - \xb'\|}.\]
 Next, due to convexity of the function $g$,
 \begin{align*}
 	g(\xb')&\ge g(\yb^*) + (\nabla g(\yb^*), \xb'-\yb^*) \\
 	& = -\rho_{t-1} + \frac{\|\yb^* - \xb'\|^2}{ \lambda \|\yb^* - \xb'\|}\ge -\rho_{t-1} + \frac{\|\yb^* - \xb'\|}{ \Lambda }.
 \end{align*}
 Thus, taking into account that $g(\xb') \le -\rho_{t}$, we obtain
 \[d =\|\yb^* - \xb'\|\le \Lambda (\rho_{t-1}-\rho_t) = O(|\rho_{t-1}-\rho_t |).\]
 \end{proof}

 With the above lemma in place, we proceed one step closer to prove Proposition~\ref{prop:zt}. Our approach for this proof is to bound $\| \zb(t+1)- \by(t) \|$ by the previous terms in the sequence, namely, $\| \zb(t) - \by(t-1)\|$ and apply a well-established lemma on convergence of a random sequence, namely,   \cite[Lemma 10, page 49]{polyak} (see Theorem~\ref{th:polyak_lem11} provided in Appendix \ref{app:supporting}). To do so though, first we need to bound the variations of $\by(t)$ as below.

 {
 \begin{lem}\label{lem:Koshal}
 	Under Assumptions~\ref{assum:convex} and~\ref{assum:Lipschitz}, the Tikhonov sequence $\by(t)$ defined in \eqref{eq:yt} satisfies
 	\begin{align*}
 		&\|\by(t)-\by(t-1)\|^2= O\left(\frac{|\e_t-\e_{t-1}|^2}{\e^{2}_t} + \frac{|\rho_t-\rho_{t-1}|^2}
 		{\e^{4+2\varepsilon}_t}\right),
 	\end{align*}
 	for any $\varepsilon>0$.
 \end{lem}}

 {
 \begin{proof}
 	Let us introduce the sequence $\hat{\by}(t) = SOL((1-\rho_{t-1}) \Ab, \Mb + \e_{t} \mathbf{I}_{Nd})$.
 	Due to the triangle inequality,
 	\begin{align}\label{eq:treq}
 		\|\by(t-1) - \by(t)\|\le\|\by(t-1)-\hat{\by}(t)\|+\|\hat{\by}(t)-\by(t)\|.
 	\end{align}
 	1) We start by analyzing $\|\by(t-1)-\hat{\by}(t)\|$.
 	Taking into account the definition of $\by(t-1)\in(1-\rho_{t-1}) \Ab$, $\hat{\by}(t)\in (1-\rho_{t-1}) \Ab$, we conclude that
 	\[\la\Mb(\by(t-1)) +\e_{t-1}\by(t-1), \hat{\by}(t) - \by(t-1)\ra\ge 0,\]
 	\[\la\Mb(\hat{\by}(t)) +\e_{t}\hat{\by}(t), \by(t-1) - \hat{\by}(t)\ra\ge 0,\]
 	By summing up two inequalities above, we obtain
 	\begin{align*}
 		0&\le  \la\Mb(\by(t-1)) +\e_{t-1}\by(t-1), \hat{\by}(t) - \by(t-1)\ra + \la\Mb(\hat{\by}(t)) +\e_{t}\hat{\by}(t), \by(t-1) - \hat{\by}(t)\ra\cr
 		&= \la\Mb(\by(t-1)) +\e_{t-1}\by(t-1), \hat{\by}(t) - \by(t-1)\ra + \la\Mb(\hat{\by}(t)) +\e_{t-1}\hat{\by}(t), \by(t-1) - \hat{\by}(t)\ra\cr
 		&\qquad - \la\Mb(\hat{\by}(t)\ra +\e_{t-1}\hat{\by}(t), \by(t-1) - \hat{\by}(t)\ra + \la\Mb(\hat{\by}(t)) +\e_{t}\hat{\by}(t), \by(t-1) - \hat{\by}(t)\ra\cr
 		&\le -\e_{t-1}\|\hat{\by}(t) - \by(t-1)\|^2 + \la(\e_{t}-\e_{t-1})\hat{\by}(t), \by(t-1) - \hat{\by}(t)\ra\cr
 		&\le -\e_{t-1}\|\hat{\by}(t) - \by(t-1)\|^2 + |\e_{t}-\e_{t-1}|\|\hat{\by}(t)\| \|\by(t-1) - \hat{\by}(t)\|,
 	\end{align*}
 	where the second inequality we used monotonicity of $\Mb$ and the last one is due to the Cauchy-Schwarz inequality.
 	Next, as $\Ab$ is compact (see Assumption~\ref{assum:convex}), we conclude existence of some finite constant $Y$ such that $\|\hat{\by}(t)\|\le Y$ for all $t$. Thus,
 	\begin{align}\label{eq:term1}
 		\|\hat{\by}(t) - \by(t-1)\|\le Y\frac{|\e_{t}-\e_{t-1}|}{\e_{t-1}}.
 	\end{align}
 	2) Next, we focus on the term $\|\hat{\by}(t)-\by(t)\|$ in~\eqref{eq:treq}. As for any $\omega_t>0$,
 	\[\hat{\by}(t) = \Proj_{(1-\rho_{t-1})\Ab}[\hat{\by}(t) - \omega_t(\Mb(\hat{\by}(t))+\e_{t}\hat{\by}(t))],\]
 	\[\by(t) = \Proj_{(1-\rho_{t})\Ab}[\by(t) - \omega_t(\Mb(\by(t))+\e_{t}\by(t))].\]
 	We consider $\omega_t$ such that $\omega_t<\e_t$.
 	\begin{align*}
 		&\|\by(t)-\hat{\by}(t)\|\cr
 		&\le\|\Proj_{(1-\rho_t)\Ab}[\by(t) - \omega_t(\Mb(\by(t))+\e_t\by(t))]-\Proj_{(1-\rho_t)\Ab}[\tilde{\by}(t)]\|\cr
 		&+ \|\Proj_{(1-\rho_t)\Ab}[\tilde{\by}(t)]-\Proj_{(1-\rho_{t-1})\Ab}[\tilde{\by}(t)]\|,
 	\end{align*}
 	where $\tilde{\by}(t)=\hat{\by}(t) - \omega_{t}(\Mb(\hat{\by}(t))+\e_{t}\hat{\by}(t))$.
 	Thus, due to  the non-expansion of the projection operator and  Lemma~\ref{lem:A_r}, we obtain
 	\begin{align*}
 		\|&\by(t)-\hat{\by}(t)\|\le\|(1-\omega_t\e_t)(\by(t)-\hat{\by}(t))-\omega_t(\Mb(\by(t))-\Mb(\hat{\by}(t)))\| \cr
 		&+O(|\rho_t-\rho_{t-1}|).
 	\end{align*}
 	We continue to estimate the term $\|(1-\omega_t\e_t)(\by(t)-\hat{\by}(t))-\omega_t(\Mb(\by(t))-\Mb(\hat{\by}(t)))\|$. Given that $\Mb$ is monotone and Lipschitz, we conclude that
 	\begin{align*}
 		&\|(1-\omega_t\e_t)(\by(t)-\hat{\by}(t))-\omega_t(\Mb(\by(t))-\Mb(\hat{\by}(t)))\|^2\cr
 		&\le(1-\omega_t\e_t)^2\|\by(t)-\hat{\by}(t)\|^2+\omega_t^2L^2\|\by(t)-\hat{\by}(t)\|^2 \cr
 		&\le (1-\omega_t\e_t)\|\by(t)-\hat{\by}(t)\|^2,
 	\end{align*}
 	where the last inequality is due to the fact that $\omega_t<\e_t$ and, thus, $\omega_t^2(\e_t^2+L^2)\le \omega_t\e_t$.
 	Hence, taking into account that $(1-\omega_t\e_t)^{1/2} \sim 1-0.5\omega_t\e_t$ for all sufficiently large $t$, we obtain
 	\begin{align*}
 		\|\by(t)-\hat{\by}(t)\|&\le(1-0.5\omega_t\e_t)\|\by(t)-\hat{\by}(t)\|+O(|\rho_t-\rho_{t-1}|),
 	\end{align*}
 	which implies that
 	\begin{align}\label{eq:term2}
 		\|&\by(t)-\hat{\by}(t)\|=O\left(\frac{|\rho_t-\rho_{t-1}|}{\omega_t\e_t}\right),
 	\end{align}
 	By taking $\omega_t = \e_t^{1+\varepsilon}$ with $\varepsilon>0$ and using the inequalities~\eqref{eq:treq},~\eqref{eq:term1},~\eqref{eq:term2}, we conclude the result.
 \end{proof}
 }

 {
 \begin{lem}\label{lem:Koshal1}
 	Assume that the interior of the set $\Ab$, $\mbox{int}(\Ab)$, is not empty and the least-norm Nash equilibria are contained in $\mbox{int}(\Ab)$. Then under Assumptions~\ref{assum:convex} and~\ref{assum:Lipschitz}, the Tikhonov sequence $\by(t)$ defined in \eqref{eq:yt} satisfies
 	\begin{align*}
 		&\|\by(t)-\by(t-1)\|^2= O\left(\frac{|\e_t-\e_{t-1}|^2}{\e^{2}_t}\right),
 	\end{align*}
 	for all sufficiently large $t$.
 \end{lem}}
 {
 \begin{proof}
 	Let $\ab^*\in\mbox{int}(\Ab)$ be the least-norm Nash equilibrium which is the limit of the sequence $\by(t)$ (see Proposition~\ref{th:Tikhonov}).
 	As $\ab^*\in\mbox{int}(\Ab)$, we conclude existence of $R>0$ such that the ball $B_{R}(\ab^*)$ with the center at $\ab^*$ and the radius $R$ is contained in the set $\Ab$, i.e. $B_{R}(\ab^*)\subseteq \Ab$. Thus,
 	\[\mbox{dist}(\ab^*,\partial \Ab)\ge R,\]
 	where $\mbox{dist}(x,\partial \Ab)$ denotes the distance between the point $x$ and the boundary $\partial \Ab$ of the compact set $\Ab$.
 	Moreover, as $\lim_{t\to\infty}\by(t)=\ab^*$, there exists such $T_R$ that for all $t\ge T_R$ we have $\|\by(t)-\ab^*\|\le \frac{R}{2}$.
 	Thus, due to the triangle inequality $\mbox{dist}(\ab^*,\partial \Ab)\le\|\by(t)-\ab^*\|+\mbox{dist}(\by(t),\partial \Ab)$,  we obtain
 	\[\mbox{dist}(\by(t),\partial \Ab)\ge\mbox{dist}(\ab^*,\partial \Ab)-\|\by(t)-\ab^*\|\ge \frac{R}{2}\]
 	for such $t\ge T_R$. Hence, for such $t$
 	\[\by(t)\in(1-R/2)\Ab\subseteq(1-\rho)\Ab\]
 	for all $\rho<\frac{R}{2}$.
 	As $\rho_t\to 0$ there exists $T_{\rho}$ such that $\rho_t<\frac{R}{2}$.
 	Next, let $T = \max(T_R,T_{\rho})$. Then for all $t>T$ and all $t'>T$ we have
 	\[\by(t) \in (1-R/2)\Ab \subseteq(1-\rho_{t'})\Ab.\]
 	In words, there exists a finite time $T$ after which the sequence $\by(t)$ is contained in the set $(1-R/2)\Ab$. This set, in its turn, is contained in the shrunk sets $(1-\rho_{t'})\Ab$ for any $t'>T$.
 	Hence, for any $t>T+1$ we conclude that $\by(t)\in (1-R/2)\Ab \subseteq(1-\rho_{t-1})\Ab$ and, thus,
 	\[(\Mb(\by(t-1)) +\e_{t-1}\by(t-1), \by(t) - \by(t-1))\ge 0,\]
 	\[(\Mb(\by(t)) +\e_{t}\by(t), \by(t-1) - \by(t))\ge 0.\]
 	By repeating the Part 1 of the proof of Lemma~\ref{lem:Koshal}, we conclude the result.
 \end{proof}
 }

 \begin{proof}(Proof of Proposition~\ref{prop:zt})
 We demonstrate that $\|\zb(t) - \by(t-1)\|\to 0$ as $t\to \infty$. Then, the result follows from Proposition~\ref{th:Tikhonov}.

 From Lemma \ref{lem:VI_sol}, since $\by(t) = \Proj_{(1-\rho_t)\Ab}[\by(t) - \g_t(\Mb(\by(t))+\e_t\by(t))]$, from definition of $\zb(t+1)$ (see~\eqref{eq:crmd}) and the non-expansion of the projection operator we obtain
 \begin{align*}
 	\|\zb(t+1) - \by(t)\|^2\le& (1-2\e_t\g_t)\|\zb(t) - \by(t)\|^2 - 2\g_t\la\zb(t) - \by(t), \Mb(\zb(t)) - \Mb(\by(t))\ra \cr
 	& + \gamma_t^2\|\Mb(\zb(t)) - \Mb(\by(t)) + \e_t(\zb(t) - \by(t))\|^2\cr
 	\le&(1-2\e_t\g_t + 2L^2\g^2_t + 2\e^2_t\g^2_t)\|\zb(t) - \by(t)\|^2\cr
 	\le& (1-\e_t\g_t)\|\zb(t) - \by(t)\|^2,
 \end{align*}
 where the last inequality is due to the fact that $\sum_{t=1}^{\infty}\e_t\g_t = \infty$, $\sum_{t=1}^{\infty}\g^2_t< \infty$ (see Assumption~\ref{assum:parameter2}), and, thus, there exists $t_0$ such that $2L^2\g^2_t + 2\e^2_t\g^2_t\le \e_t\g_t$ for  $t>t_0$.
 Hence, for sufficiently large $t$
 \begin{align*}
 	\nonumber
 	\|&\zb(t+1) - \by(t)\|^2\le (1+0.5\e_t\g_t)(1-\e_t\g_t)\|\zb(t) - \by(t-1)\|^2 \\
 	\nonumber
 	&+ (1+\frac{2}{\e_t\g_t})\|\by(t) - \by(t-1)\|^2\\
 	\nonumber
 	\le&(1-0.5\e_t\g_t)\|\zb(t) - \by(t-1)\|^2 \\
 	\nonumber
 	&+ O\left(\frac{|\e_t-\e_{t-1}|^2}{\e^3_t\g_t}+\frac{|\rho_t-\rho_{t-1}|^2}{\e^{5+2\varepsilon}_t\g_t}\right),
 \end{align*}
 where in the first inequality we used the fact that $\|a+b\|^2\le (1+\theta)\|a\|^2 + \left(1+\frac{1}{\theta}\right)\|b\|^2$ for any two vectors $a$ and $b$ and any constant $\theta>0$ and in last inequality we used Lemma~\ref{lem:Koshal}. Next, by taking into account Assumption~\ref{assum:parameter2} and the result in Theorem~\ref{th:polyak_lem11} (see \ref{app:supporting}) for the sequence $\|\zb(t) - \by(t-1)\|$, we conclude that $\zb(t) \to \by(t-1)$ as $t\to \infty$. Thus, according to Proposition~\ref{th:Tikhonov}, the sequence $\zb(t)$ converges to the least-norm solution of $VI(\Ab,\Mb)$.
 {Note that in the case when the least-norm solution $\ab^*$ of $VI(\Ab,\Mb)$ is such that $\ab^*\in\mbox{int}(\Ab)$, Lemma~\ref{lem:Koshal1} is applicable, which implies
 	\begin{align*}
 		\nonumber
 		\|&\zb(t+1) - \by(t)\|^2\le (1+0.5\e_t\g_t)(1-\e_t\g_t)\|\zb(t) - \by(t-1)\|^2 \\
 		\le&(1-0.5\e_t\g_t)\|\zb(t) - \by(t-1)\|^2 + O\left(\frac{|\e_t-\e_{t-1}|^2}{\e^3_t\g_t}\right),
 	\end{align*}
 	and, thus the condition~\ref{itm:omega} in Assumption~\ref{assum:parameter2}  can be replaced by~\ref{itm:omega1}.}
 \end{proof}
 
\section{Proof of Lemma \ref{lem:muRegret}}\label{app:muRegret}
 \begin{proof}
 \[R^i(T,\ab^i) - R_{\bmu}^i(T,\ab^i) = \sum_{t=1}^{T}[J_t^i(\ab^i(t)) - J_{t,\bmu}^i(\bmu^i(t))] + \sum_{t=1}^{T}[J_{t,\bmu}^i(\ab^i) - J_t^i(\ab^i)].\]
 We start by estimating the expectation of the terms $J_t^i(\ab^i(t)) - J_{t,\bmu}^i(\bmu^i(t))$ in the sum above.
 \begin{align*}
 	J_t^i(\ab^i(t)) - J_{t,\bmu}^i(\bmu^i(t)) = J^i(\ab(t)) - J^i(\bmu(t)) = J^i(\bxi(t)) - J^i(\bmu(t)) + J^i(\ab(t)) - J^i(\bxi(t)),
 \end{align*}
 where, as it has been defined in Section~\ref{sec:analysis}, $\bxi(t)$ is the randomly distributed vector with the independent coordinates sampled from $\EuScript N(\bmu^i(t),\sigma_t)$ and $\ab(t) = \Proj_{\Ab} \bxi(t).$ Thus,
 \begin{align}\label{eq:cond_expFull}
 	\E\{J_t^i(\ab^i(t)) - J_{t,\bmu}^i(\ab^i(t))|\EuScript F_t\} = \E\{J^i(\bxi(t)) - J^i(\bmu(t))|\EuScript F_t\} + \E\{J^i(\ab(t)) - J^i(\bxi(t))|\EuScript F_t\}.
 \end{align}

 We proceed with estimating $\E\{J^i(\bxi(t)) - J^i(\bmu(t))|\EuScript F_t\}$:
 \begin{align}\label{eq:diffJtmu}
 	&\E\{J^i(\bxi(t)) - J^i(\bmu(t))|\EuScript F_t\} = \int_{\R^{Nd}}[J^i(\bx)-J^i(\bmu)] p(\bx;\bmu, \sigma_t)d\bx \cr
 	&\le\int_{\R^{Nd}}(\nabla J^i(\bmu),\bx - \bmu) p(\bx;\bmu, \sigma_t)d\bx +\int_{\R^{Nd}}\|\nabla^2 {J^i}(\tilde{\bmu})\| \|\bx-\bmu\|^2  p(\bx;\bmu, \sigma_t) d\bx\cr
 	&= \int_{\R^{Nd}}\|\nabla^2 {J^i}(\tilde{\bmu})\| \|\bx-\bmu\|^2  p(\bx;\bmu, \sigma_t) d\bx,
 \end{align}
 where $p = \prod_{i=1}^N p^i$ and $p^i$ is defined in~\eqref{eq:density} and $\tilde{\bmu} = \bmu + \theta (\bx-\bmu)$, for some $\theta\in[0,1]$ given $\bx$ and $\bmu$, and $\nabla^2 {J^i}(\tilde{\bmu})$ is the Hessian matrix of the function $J^i$ at the point (which exists according to Assumption~\ref{assum:Lipschitz}~\ref{itm:lip}).
 The first inequality above is obtained by the mean value theorem (which is applicable again according to Assumption~\ref{assum:Lipschitz}~\ref{itm:lip}) and the second one is due to the definition of the density function $p$ implying
 \[\int_{\R^{Nd}}(\nabla J^i(\bmu),\bx - \bmu) p(\bx;\bmu, \sigma_t)d\bx = 0.\]
 By applying the H\"older's inequality and the boundedness of the integral $\int_{\R^{Nd}}\|\nabla^2 {J^i}(\tilde{\bmu})\| \|\bx-\bmu\|^2  p(\bx;\bmu, \sigma_t) d\bx$ from Assumption \ref{assum:Lipschitz} Part \ref{itm:growth} (see~\eqref{eq:boundR} for an analogous result), we obtain
 \begin{align}\label{eq:cond_expTerm1}
 	&\E\{J^i(\bxi(t)) - J^i(\bmu(t))|\EuScript F_t\} = O(\sigma_t^2),
 \end{align}
 Next, to address  $\E\{J^i(\ab(t)) - J^i(\bxi(t))|\EuScript F_t\}$ in~\eqref{eq:cond_expFull}, repeating the proof of Lemma~\ref{lem:projTerm} we obtain
 \begin{align}\label{eq:cond_expTerm2}
 	&\E\{J^i(\ab(t)) - J^i(\bxi(t))|\EuScript F_t\}  =O\left(\left(\frac{e^{-\frac{{\rho}_t^2}{2\sigma_t^2}}}{\sigma_t^{Nd}}\right)^{\frac{1}{2}}\right).
 \end{align}
 Thus, bringing the inequalities~\eqref{eq:cond_expFull}, \eqref{eq:cond_expTerm1}, \eqref{eq:cond_expTerm2} together, and taking the full expectation, we conclude the result due to the condition  $\sigma_t = \frac{1}{t^s}$, $\rho_t = \frac{1}{t^r}$,  with $0< g,e,s,r$, $r<s$.
 \end{proof}

 \section{Proof of Lemma \ref{lem:mixedstr_cost}}\label{app:mixed_cost}
 \begin{proof}
  Part ~\ref{itm:convex_smooth}:  To show convexity of $\tilde J^i_{t,\bmu}$ consider $\bmu^i_1, \bmu^i_2 \in \Ab^i$, $a_1, a_2 \in \R$. From the definition of $\tilde J^i_{t,\bmu}$ (see~\eqref{eq:tildeJ_t}) and definition of $p^j$ in~\eqref{eq:density}, we have
 \begin{align}\label{eq:eq1}
 	&\tilde J^i_{t,\bmu}(a\bmu^i_1 + (1-a)\bmu^i_2) \cr
 	&= \frac{1}{(2\pi\sigma^2)^{d/2}}\int_{\R^{Nd}}J^i(\bx) \exp\left\{-\frac{\|\bx^i-a\bmu^i_1 - (1-a)\bmu^i_2\|^2}{2\sigma^2}\right\}\prod_{j\ne i}p^j(\bx^j;\bmu^j(t),\sigma_t)d\bx.
 \end{align}
 Substituting $\by^i = \bx^i-a\bmu^i_1 - (1-a)\bmu^i_2$, $\by^{-i} = \bx^{-i}$, we get
 \begin{align}\label{eq:eq2}
 	&\int_{\R^{Nd}}J^i(\bx) \exp\left\{-\frac{\|\bx^i-a\bmu^i_1 - (1-a)\bmu^i_2\|^2}{2\sigma_t^2}\right\}\prod_{j\ne i}p^j(\bx^j;\bmu^j(t),\sigma_t)d\bx\cr
 	& = \int_{\R^{Nd}}J^i(\by^i + a\bmu^i_1 +(1-a)\bmu^i_2, \by^{-i}) \exp\left\{-\frac{\|\by^i\|^2}{2\sigma_t^2}\right\}\prod_{j\ne i}p^j(\by^j;\bmu^j(t),\sigma_t)d\by\cr
 	& = \int_{\R^{Nd}}J^i(a(\by^i + \bmu^i_1) +(1-a)(\by^i + \bmu^i_2)) \exp\left\{-\frac{\|\by^i\|^2}{2\sigma_t^2}\right\}\prod_{j\ne i}p^j(\by^j;\bmu^j(t),\sigma_t)d\by\cr
 	& \le a\int_{\R^{Nd}}J^i(\by^i + \bmu^i_1,\by^{-i})\exp\left\{-\frac{\|\by^i\|^2}{2\sigma_t^2}\right\}\prod_{j\ne i}p^j(\by^j;\bmu^j(t),\sigma_t)d\by\cr
 	&\quad +(1-a)\int_{\R^{Nd}}J^i(\by^i + \bmu^i_2,\by^{-i}) \exp\left\{-\frac{\|\by^i\|^2}{2\sigma_t^2}\right\}\prod_{j\ne i}p^j(\by^j;\bmu^j(t),\sigma_t)d\by\cr
 	& = a\int_{\R^{Nd}}J^i_{t}(\bx)\exp\left\{-\frac{\|\bx - \bmu^i_1\|^2}{2\sigma_t^2}\right\}\prod_{j\ne i}p^j(\bx^j;\bmu^j(t),\sigma_t)d\bx\cr
 	& \quad +(1-a)\int_{\R^{Nd}}J^i(\bx) \exp\left\{-\frac{\|\bx- \bmu^i_2\|^2}{2\sigma_t^2}\right\}\prod_{j\ne i}p^j(\bx^j;\bmu^j(t),\sigma_t)d\bx
 \end{align}
 where the first inequality is due to $J^i(\bx)$ being convex in $\bx^i$ by Assumption~\ref{assum:convex}. Hence, by combining~\eqref{eq:eq1} and~\eqref{eq:eq2}, we conclude
 \begin{align*}
 	\tilde J^i_{t,\bmu}(a\bmu^i_1 + (1-a)\bmu^i_2)  \leq a\tilde J^i_{t,\bmu}(\bmu^i_1) + (1-a)\tilde J^i_{t,\bmu}(\bmu^i_2).
 \end{align*}
 The fact that for any fixed $t$ the gradient $\nabla\tilde J^i_{t,\bmu}(\cdot)$ is bounded over the compact set $\Ab$ follows directly from its continuity over $\Ab$ implied by differentiability of the function $\tilde J^i_{t,\bmu}(\cdot)$ (see Lemma ~\ref{lem:sample_grad} Part~\ref{itm:diff_int} and \cite[Chapter 17]{zorich}).

 Part \ref{itm:bound_smooth}: According to the definitions of the functions $J^i_{t,\bmu}$ and $\tilde J^i_{t,\bmu}$ (see~\eqref{eq:Jtmu}, \eqref{eq:tildeJ_t})
 \begin{align*}
 	\; &J^i_{t,\bmu}(\bmu^i) - \tilde J^i_{t,\bmu}(\bmu^i)= \int_{\R^{Nd}}[J^i(\bmu^i,\bmu^{-i}(t))-J^i(\bx)] p^i(\bx^i;\bmu^i, \sigma_t)\prod_{j\ne i}p^j(\bx^j;\bmu^j(t), \sigma_t)d\bx.
 \end{align*}
 Thus, analogously to~\eqref{eq:diffJtmu}, we conclude that  $J^i_{t,\bmu}(\bmu^i) - \tilde J^i_{t,\bmu}(\bmu^i) = O(\sigma_t).$
 \end{proof}

\section{Supporting Theorems}\label{app:supporting}
The following result related to the convergence of the stochastic process is proven in Lemma 10 (page 49) in \cite{polyak}.
\begin{theorem}\label{th:polyak_lem11}
Let $v_0, \ldots, v_t$ be a sequence of random variables, $v_t\ge 0$, $\E v_0<\infty$ and let
\[\E \{v_{t+1}| \EuScript F_t\} \le (1-\alpha_t)v_t + \phi_t,\]
where $\EuScript F_t$ is the $\sigma$-algebra generated by the random variables $\{v_0,\ldots,v_t\}$, $0<\alpha_t<1$, $\sum_{t=0}^{\infty}\alpha_t = \infty$,
$\beta_t\ge 0$, $\sum_{t=0}^{\infty}\beta_t < \infty$. 
Then $v_t\to 0$ almost surely, $\E v_t\to 0$ as $t\to \infty$.
\end{theorem}

\end{document}